\documentclass[12pt]{article}
\usepackage{latexsym, amssymb}

\DeclareMathAlphabet\gothic{U}{euf}{m}{n}

\makeatletter
\def\eqnarray{\stepcounter{equation}\let\@currentlabel=\theequation
\global\@eqnswtrue
\tabskip\@centering\let\\=\@eqncr
$$\halign to \displaywidth\bgroup\hfil\global\@eqcnt\z@
  $\displaystyle\tabskip\z@{##}$&\global\@eqcnt\@ne
  \hfil$\displaystyle{{}##{}}$\hfil
  &\global\@eqcnt\tw@ $\displaystyle{##}$\hfil
  \tabskip\@centering&\llap{##}\tabskip\z@\cr}

\def\endeqnarray{\@@eqncr\egroup
      \global\advance\c@equation\m@ne$$\global\@ignoretrue}

\def\@yeqncr{\@ifnextchar [{\@xeqncr}{\@xeqncr[5pt]}}
\makeatother

\begin{document}
\bibliographystyle{tom}

\newtheorem{lemma}{Lemma}[section]
\newtheorem{thm}[lemma]{Theorem}
\newtheorem{cor}[lemma]{Corollary}
\newtheorem{voorb}[lemma]{Example}
\newtheorem{rem}[lemma]{Remark}
\newtheorem{prop}[lemma]{Proposition}
\newtheorem{stat}[lemma]{{\hspace{-5pt}}}
\newtheorem{obs}[lemma]{Observation}
\newtheorem{defin}[lemma]{Definition}

\newenvironment{remarkn}{\begin{rem} \rm}{\end{rem}}
\newenvironment{exam}{\begin{voorb} \rm}{\end{voorb}}
\newenvironment{defn}{\begin{defin} \rm}{\end{defin}}
\newenvironment{obsn}{\begin{obs} \rm}{\end{obs}}

\newenvironment{emphit}{\begin{itemize} }{\end{itemize}}

\newcommand{\gota}{\gothic{a}}
\newcommand{\gotb}{\gothic{b}}
\newcommand{\gotc}{\gothic{c}}
\newcommand{\gote}{\gothic{e}}
\newcommand{\gotf}{\gothic{f}}
\newcommand{\gotg}{\gothic{g}}
\newcommand{\gothh}{\gothic{h}}
\newcommand{\gotk}{\gothic{k}}
\newcommand{\gotm}{\gothic{m}}
\newcommand{\gotn}{\gothic{n}}
\newcommand{\gotp}{\gothic{p}}
\newcommand{\gotq}{\gothic{q}}
\newcommand{\gotr}{\gothic{r}}
\newcommand{\gots}{\gothic{s}}
\newcommand{\gotu}{\gothic{u}}
\newcommand{\gotv}{\gothic{v}}
\newcommand{\gotw}{\gothic{w}}
\newcommand{\gotz}{\gothic{z}}
\newcommand{\gotA}{\gothic{A}}
\newcommand{\gotB}{\gothic{B}}
\newcommand{\gotG}{\gothic{G}}
\newcommand{\gotL}{\gothic{L}}
\newcommand{\gotS}{\gothic{S}}
\newcommand{\gotT}{\gothic{T}}

\newcommand{\mn}{\marginpar{\hspace{1cm}*} }
\newcommand{\mnn}{\marginpar{\hspace{1cm}**} }

\newcommand{\mnq}{\marginpar{\hspace{1cm}*???} }
\newcommand{\mnnq}{\marginpar{\hspace{1cm}**???} }

\newcounter{teller}
\renewcommand{\theteller}{\Roman{teller}}
\newenvironment{tabel}{\begin{list}%
{\rm \bf \Roman{teller}.\hfill}{\usecounter{teller} \leftmargin=1.1cm
\labelwidth=1.1cm \labelsep=0cm \parsep=0cm}
                      }{\end{list}}

\newcounter{tellerr}
\renewcommand{\thetellerr}{(\roman{tellerr})}
\newenvironment{subtabel}{\begin{list}%
{\rm  (\roman{tellerr})\hfill}{\usecounter{tellerr} \leftmargin=1.1cm
\labelwidth=1.1cm \labelsep=0cm \parsep=0cm}
                         }{\end{list}}
\newenvironment{ssubtabel}{\begin{list}%
{\rm  (\roman{tellerr})\hfill}{\usecounter{tellerr} \leftmargin=1.1cm
\labelwidth=1.1cm \labelsep=0cm \parsep=0cm \topsep=1.5mm}
                         }{\end{list}}

\newcommand{\Ni}{{\bf N}}
\newcommand{\Ri}{{\bf R}}
\newcommand{\Ci}{{\bf C}}
\newcommand{\Ti}{{\bf T}}
\newcommand{\Zi}{{\bf Z}}
\newcommand{\Fi}{{\bf F}}

\newcommand{\proof}{\mbox{\bf Proof} \hspace{5pt}} 
\newcommand{\remark}{\mbox{\bf Remark} \hspace{5pt}}
\newcommand{\ruimte}{\vskip10.0pt plus 4.0pt minus 6.0pt}

\newcommand{\simh}{{\stackrel{{\rm cap}}{\sim}}}
\newcommand{\ad}{{\mathop{\rm ad}}}
\newcommand{\Ad}{{\mathop{\rm Ad}}}
\newcommand{\Aut}{\mathop{\rm Aut}}
\newcommand{\arccot}{\mathop{\rm arccot}}
\newcommand{\capp}{{\mathop{\rm cap}}}
\newcommand{\rcapp}{{\mathop{\rm rcap}}}
\newcommand{\diam}{\mathop{\rm diam}}
\newcommand{\divv}{\mathop{\rm div}}
\newcommand{\codim}{\mathop{\rm codim}}
\newcommand{\RRe}{\mathop{\rm Re}}
\newcommand{\IIm}{\mathop{\rm Im}}
\newcommand{\Tr}{{\mathop{\rm Tr}}}
\newcommand{\Vol}{{\mathop{\rm Vol}}}
\newcommand{\card}{{\mathop{\rm card}}}
\newcommand{\supp}{\mathop{\rm supp}}
\newcommand{\sgn}{\mathop{\rm sgn}}
\newcommand{\essinf}{\mathop{\rm ess\,inf}}
\newcommand{\esssup}{\mathop{\rm ess\,sup}}
\newcommand{\Int}{\mathop{\rm Int}}
\newcommand{\Leibniz}{\mathop{\rm Leibniz}}
\newcommand{\lcm}{\mathop{\rm lcm}}
\newcommand{\loc}{{\rm loc}}

\newcommand{\mod}{\mathop{\rm mod}}
\newcommand{\spann}{\mathop{\rm span}}
\newcommand{\one}{1\hspace{-4.5pt}1}

\newcommand{\DWR}{}

\hyphenation{groups}
\hyphenation{unitary}

\newcommand{\tfrac}[2]{{\textstyle \frac{#1}{#2}}}

\newcommand{\cb}{{\cal B}}
\newcommand{\cc}{{\cal C}}
\newcommand{\cd}{{\cal D}}
\newcommand{\ce}{{\cal E}}
\newcommand{\cf}{{\cal F}}
\newcommand{\ch}{{\cal H}}
\newcommand{\ci}{{\cal I}}
\newcommand{\ck}{{\cal K}}
\newcommand{\cl}{{\cal L}}
\newcommand{\cm}{{\cal M}}
\newcommand{\cn}{{\cal N}}
\newcommand{\co}{{\cal O}}
\newcommand{\cs}{{\cal S}}
\newcommand{\ct}{{\cal T}}
\newcommand{\cx}{{\cal X}}
\newcommand{\cy}{{\cal Y}}
\newcommand{\cz}{{\cal Z}}

\newcommand{\wtozp}{W^{1,2}\raisebox{10pt}[0pt][0pt]{\makebox[0pt]{\hspace{-34pt}$\scriptstyle\circ$}}}
\newlength{\hightcharacter}
\newlength{\widthcharacter}
\newcommand{\covsup}[1]{\settowidth{\widthcharacter}{$#1$}\addtolength{\widthcharacter}{-0.15em}\settoheight{\hightcharacter}{$#1$}\addtolength{\hightcharacter}{0.1ex}#1\raisebox{\hightcharacter}[0pt][0pt]{\makebox[0pt]{\hspace{-\widthcharacter}$\scriptstyle\circ$}}}
\newcommand{\cov}[1]{\settowidth{\widthcharacter}{$#1$}\addtolength{\widthcharacter}{-0.15em}\settoheight{\hightcharacter}{$#1$}\addtolength{\hightcharacter}{0.1ex}#1\raisebox{\hightcharacter}{\makebox[0pt]{\hspace{-\widthcharacter}$\scriptstyle\circ$}}}
\newcommand{\scov}[1]{\settowidth{\widthcharacter}{$#1$}\addtolength{\widthcharacter}{-0.15em}\settoheight{\hightcharacter}{$#1$}\addtolength{\hightcharacter}{0.1ex}#1\raisebox{0.7\hightcharacter}{\makebox[0pt]{\hspace{-\widthcharacter}$\scriptstyle\circ$}}}

 \thispagestyle{empty}

\vspace*{1.0cm}
\begin{center}
{\Large{\bf Markov uniqueness }}\\[3mm] 
{\Large{\bf of degenerate elliptic operators  }}  \\[5mm]
\large  Derek W. Robinson$^1$ and Adam Sikora$^2$\\[2mm]

\normalsize{November 2009}
\end{center}

\vspace{5mm}

\begin{center}
{\bf Abstract}
\end{center}

\begin{list}{}{\leftmargin=1.8cm \rightmargin=1.8cm \listparindent=10mm 
   \parsep=0pt}
   \item
   Let $\Omega$ be an open subset of $\Ri^d$ and $H_\Omega=-\sum^d_{i,j=1}\partial_i\,c_{ij}\,\partial_j$
a   second-order partial differential operator on $L_2(\Omega)$ with domain $C_c^\infty(\Omega)$  
where  the coefficients $c_{ij}\in W^{1,\infty}(\Omega)$ are real symmetric
and  $C=(c_{ij})$ is a strictly positive-definite matrix over $\Omega$.
 In particular, $H_\Omega$ is locally strongly elliptic.

We analyze the submarkovian extensions of $H_\Omega$, i.e.\ the self-adjoint extensions which 
generate submarkovian semigroups.
Our main result establishes  that $H_\Omega$ is Markov unique, i.e.\ it has a unique submarkovian extension,
  if and only if   $\capp_\Omega(\partial\Omega)=0$ where 
$\capp_\Omega(\partial\Omega)$ is the capacity of the  boundary of $\Omega$ measured with respect to $H_\Omega$.
The second main result establishes that Markov uniqueness of $H_\Omega$ is equivalent to the  semigroup generated
by the Friedrichs extension of $H_\Omega$ being conservative.

\end{list}

\vfill

\noindent AMS Subject Classification: 47B25, 47D07, 35J70.

\vspace{0.5cm}

\noindent
\begin{tabular}{@{}cl@{\hspace{10mm}}cl}
1. & Centre for Mathematics & 
  2. & Department of Mathematics  \\
&\hspace{15mm} and its Applications  & 
  &Macquarie University  \\
& Mathematical Sciences Institute & 
  & Sydney, NSW 2109  \\
& Australian National University& 
  & Australia \\
& Canberra, ACT 0200  & {}
  & \\
& Australia & {}
  & \\
& derek.robinson@anu.edu.au & {}
  &sikora@ics.mq.edu.au \\
\end{tabular}

\newpage

\setcounter{page}{1}

\section{Introduction}\label{S1}

The Markov uniqueness problem \cite{Ebe}  consists of finding  conditions which ensure that  a diffusion operator has  a unique
submarkovian extension, i.e.\ an extension that generates a submarkovian semigroup.
An operator with this property is said to be Markov unique.
Our aim is to analyze this problem for the class of   second-order, divergence-form,  elliptic operators
 with real Lipschitz continuous coefficients acting on an open subset of $\Omega$ of $\Ri^d$.
 Each of these operators has at least one submarkovian extension, the Friedrichs extension  \cite{Friedr2}.
This  extension corresponds to Dirichlet boundary conditions on $\partial\Omega$ and alternative boundary conditions can lead to  different submarkovian extensions.
 Our principal results establish that Markov uniqueness is equivalent to the boundary $\partial\Omega$  having zero capacity, Theorem~\ref{tsm1.1},
 or to conservation of probability, Theorem~\ref{tsm1.2}.

 Define $H_\Omega$ as  the positive symmetric operator  on $L_2(\Omega)$ with
 domain $D(H_\Omega)=C_c^\infty(\Omega)$ and action
\begin{equation}
H_\Omega\varphi= -\sum^d_{i,j=1}\partial_i\,c_{ij}\,\partial_j\varphi= -\sum^d_{i,j=1}c_{ij}\,\partial_i\,\partial_j\varphi-\sum^d_{i,j=1}(\partial_ic_{ij})\,\partial_j\varphi
\label{ese1.0}
\end{equation}
where the $c_{ij}=c_{ji}\in W^{1,\infty}(\Omega)$ are real, $C=(c_{ij})$
is a  non-zero, positive-definite matrix over $\Omega$ and $\partial_i=\partial/\partial x_i$.
We assume throughout that $C(x)=(c_{ij}(x))>0$ for all $x\in\Omega$.
This ensures that $H_\Omega$ is  locally strongly elliptic, i.e.\ 
 for each  compact subset $K$ of $\Omega$
there is a $\mu_K>0$ such that $C(x)\geq \mu_K I$ for all $x\in K$.
This ellipticity property is fundamental as it ensures that the various possible self-adjoint extensions of $H_\Omega$ differ only in their boundary behaviour 
(see Section~\ref{S3}).

The Markov uniqueness problem has been considered in a variety of contexts (see \cite{Ebe} for background material and an extensive survey).
It is related to a number of other uniqueness problems.
For example, the operator $H_\Omega$, which can be viewed as an operator on $L_p(\Omega)$ for each $p\in[1,\infty]$, is defined to be $L_p$-unique
if it has a unique extension which generates an $L_p$-continuous semigroup.
In particular $H_\Omega$ is $L_2$-unique if and only if it is essentially self-adjoint (see \cite{Ebe}, Corollary~1.2).
Then the self-adjoint closure is automatically submarkovian and  $H_\Omega$ is Markov unique.
Moreover, if $H_\Omega$ is $L_1$-unique then it is Markov unique (\cite{Ebe},  Lemma~1.6).
In Theorem~\ref{tsm1.2} we will establish  a converse to this statement for the class of operators under consideration.
As a byproduct of our analysis of Markov uniqueness we also  derive criteria for various other forms of uniqueness.

In the sequel we  use extensively the theory of positive closed quadratic forms and positive self-adjoint operators  (see \cite{Kat1}, Chapter~6) and the corresponding theory
of Dirichlet forms and submarkovian operators (see \cite{BH} \cite{MR} \cite{FOT}).
First we introduce the quadratic form $h_\Omega$  associated with $H_\Omega$ by 
\begin{equation}
h_\Omega(\varphi)
=\sum^d_{i,j=1} \int_\Omega dx\, c_{ij}(x)\,(\partial_i\varphi)(x)(\partial_j\varphi)(x)
\label{ese1.00}
\end{equation}
with domain $D(h_\Omega)=D(H_\Omega)=C_c^\infty(\Omega)$.
The form $h_\Omega$ is closable with respect to the graph norm $\varphi\mapsto\|\varphi\|_{D(h_\Omega)}=(h_\Omega(\varphi)+\|\varphi\|_2^2)^{1/2}$
and its closure $\overline h_\Omega$ is a Dirichlet form.
The  positive self-adjoint operator corresponding to $\overline h_\Omega$ is the Friedrichs extension $H_\Omega^F$  of $H_\Omega$.
It is automatically submarkovian. 
Moreover, it  is the largest positive self-adjoint extension of $H_\Omega$ with respect to the  usual ordering of self-adjoint operators.
Krein  \cite{Kre1}   established that $H_\Omega$ also has a smallest positive self-adjoint extension.
But the  Krein extension is not always  submarkovian.
For example, if $\Omega$ is bounded the Krein extension of the Laplacian restricted to $C_c^\infty(\Omega)$ is not submarkovian (see \cite{FOT}, Theorem~3.3.3).
Our first aim is to establish that $H_\Omega$ also has a smallest submarkovian extension.
Then the Markov uniqueness problem is reduced to finding conditions which ensure that this latter extension coincides with the Friedrichs extension
(see \cite{Ebe}, Chapter~3).

Define $l_\Omega$ by setting 
\begin{equation} l_\Omega(\varphi)=\int_\Omega dx\, \sum^d_{i,j=1} c_{ij}(x)\,(\partial_i\varphi)(x)(\partial_j\varphi)(x)
\label{esm1.1} \end{equation}
where the $\partial_i \varphi$ denote the distributional derivatives and the domain $D(l_\Omega)$ of the form is defined to be the space
of all $\varphi\in L_2(\Omega)$ for which the integral is finite.
It is clear that $l_\Omega$ is an extension of $\overline h_\Omega$ but it is not immediately obvious that $l_\Omega$ is closed 
and that the corresponding operator $L_\Omega$ is an extension of $H_\Omega$.
These properties were established in  \cite{FOT} for operators of the form (\ref{ese1.0}) but with smooth coefficients.
Our first result is a generalization for
operators with Lipschitz coefficients which also incorporates some regularity and domination properties.

Recall that the positive semigroup $S_t$ is defined to dominate the positive semigroup $T_t$ if
$S_t\varphi\geq T_t\varphi$ for all positive $\varphi\in L_2(\Omega)$ and all $t>0$.
Moreover, $D(k_\Omega)$ is defined to be an order ideal of $D(l_\Omega)$ 
if the conditions $0\leq \varphi\leq \psi$, $\psi\in D(k_\Omega)$ and $\varphi\in D(l_\Omega)$ imply $\varphi\in D(k_\Omega)$.
(See \cite{Ouh5}, Chapter~2, for these and related concepts.)

\begin{thm}\label{tsm1.0}
Let $\Omega$ be an open subset of $\Ri^d$ and $H_\Omega=-\sum^d_{i,j=1}\partial_i\,c_{ij}\,\partial_j$
a   second-order partial differential operator on $L_2(\Omega)$ with domain $C_c^\infty(\Omega)$  
where   the $c_{ij}\in W^{1,\infty}(\Omega)$ are real symmetric
and  $C(x)=(c_{ij}(x))>0$ for all $x\in \Omega$.

Then the following are true.
  \begin{tabel}
  \item\label{tsm1.0-0}
 $l_\Omega$ is a Dirichlet form and $D(l_\Omega)\cap C^\infty(\Omega)$ is a core of $l_\Omega$.
  \item\label{tsm1.0-1}
 The submarkovian  operator $L_\Omega$ associated with $l_\Omega$ is an extension of $H_\Omega$.
  \item\label{tsm1.0-2}
If  $K_\Omega$ is  any  submarkovian extension of $H_\Omega$  and $k_\Omega$   the corresponding Dirichlet form
then $l_\Omega\supseteq k_\Omega\supseteq h_\Omega$.
Therefore $0\leq L_\Omega\leq K_\Omega\leq H_\Omega^F$ in the sense of operator order.
 \item\label{tsm1.0-3}
 If $K_\Omega$ is any self-adjoint extension of $H_\Omega$ then 
 $C_c^\infty(\Omega)D(K_\Omega)\subseteq D(\overline H_\Omega)$
 and if $K_\Omega$ is a submarkovian extension then
 $C_c^\infty(\Omega)D(k_\Omega)\subseteq D(\overline h_\Omega)$.
 \item\label{tsm1.0-5}
If $K_\Omega$ is a submarkovian extension of $H_\Omega$ then  $D(\overline  h_\Omega)$ is an order ideal of $D(k_\Omega)$
and  $e^{-tK_\Omega}$ dominates $ e^{-tH_\Omega^F}$.
Moreover, $e^{-tL_\Omega}$ dominates $ e^{-tK_\Omega}$
if and only if $D(k_\Omega)$ is an order ideal of $D(l_\Omega)$.
\end{tabel}
\end{thm}

The first three statements are a generalization of Lemma~3.3.3 and Theorem~3.3.1 in \cite{FOT}.
They establish that  the operator  $L_\Omega$  is the smallest submarkovian extension of $H_\Omega$,
but not necessarily the smallest self-adjoint extension.
The fourth statement is an interior regularity property.
It establishes, in particular, that every submarkovian extension of $H_\Omega$ is a Silverstein extension in the terminology
of \cite{Take} (see \cite{Ebe},  Definition~1.4).

The third   statement of the theorem implies that $H_\Omega$ is Markov unique if and only if $l_\Omega=\overline h_\Omega$, i.e.\ if and only if $D(l_\Omega)=D(\overline h_\Omega)$.
It is this criterion that has been used extensively in the analysis of the Markov uniqueness problem (see \cite{AKR} \cite{FOT}  and  \cite{Ebe}, Chapter~3).
But the fourth statement implies that all the Dirichlet form extensions coincide in the interior of $\Omega$ and consequently differ only on the boundary.
Our first criterion for 
Markov uniqueness is   in terms of the capacity of the boundary.

The (relative) capacity   of the measurable subset $A\subset \overline \Omega$ is defined by
\begin{eqnarray}
\capp_\Omega(A)=\inf\Big\{\|\psi\|_{D(l_\Omega)}^2&&: \psi\in D(l_\Omega) \mbox{ and  there exists   an open set  }\nonumber\\[-5pt]
&& U\subset \Ri^d
\mbox{ such that } U\supseteq A
\mbox{ and } \psi=1 \mbox{ a.\ e.\ on } U\cap\Omega\Big\}\label{ecap}
\;.
\end{eqnarray}
Thus $\capp_\Omega$ is directly related to the capacity occurring in the theory of 
Dirichlet forms \cite{BH} \cite{FOT}.  

\begin{thm}\label{tsm1.1}
Under the assumptions of Theorem~$\ref{tsm1.0}$, the following conditions are equivalent:\vspace{-2mm}
\begin{tabel}
\item\label{tsm1.1-1}
$H_\Omega$ is Markov unique, 
\item\label{tsm1.1-2}
$\capp_\Omega(\partial\Omega)=0$.
\end{tabel}
\end{thm}

It should be emphasized  that there is no comparable geometric or potential-theoretic characterization of essential self-adjointness, i.e.\ $L_2$-uniqueness.
Folklore would suggest that $H_\Omega$ is $L_2$-unique if and only if the Riemannian distance to the boundary $\partial\Omega$, measured with respect to the metric $C^{-1}$, is infinite.
But this is not true in one-dimension (see Example~\ref{exsm5.1}).

Our second result on Markov uniqueness is based on a conservation property.
The submarkovian semigroup $S^F_t$ generated by the Friedrichs extension $H_\Omega^F$ is defined to be conservative on $L_\infty(\Omega)$
 if $S^F_t\one_\Omega=\one_\Omega$  for all $t\geq0$.
 
 \begin{thm}\label{tsm1.2} 
Adopt the assumptions of Theorem~$\ref{tsm1.0}$.
Let $S^F_t$ denote the semigroup generated by the Friedrichs extension $H_\Omega^F$ of $H_\Omega$.\vspace{1mm}

The following conditions are equivalent:\vspace{-1mm}
\begin{tabel}
\item\label{tsm1.2-3}
$H_\Omega$ is Markov unique,
\item\label{tsm1.2-0}
$S^F_t$ is conservative,
\item\label{tsm1.2-2}
$H_\Omega$ is $L_1$-unique.
\end{tabel}
\end{thm}

The implications \ref{tsm1.2-0}$\Leftrightarrow$\ref{tsm1.2-2}$\Rightarrow$\ref{tsm1.2-3} are already known under slightly different hypotheses.
The equivalence of Conditions~\ref{tsm1.2-0} and \ref{tsm1.2-2} was established by Davies, \cite{Dav14} Theorem~2.2, for a 
different class of second-order operators with smooth coefficients.
His proof is based on an earlier result of Azencott \cite{Aze}.
The implication  \ref{tsm1.2-2}$\Rightarrow$\ref{tsm1.2-3}  is quite general  and is given by  Lemma~1.6 in 
\cite{Ebe}.
Moreover, the implication \ref{tsm1.2-3}$\Rightarrow$\ref{tsm1.2-0} follows from \cite{Ebe}, Corollary~3.4, if $|\Omega|<\infty$.
The proof of this implication for general $\Omega$ is considerably more complicated (see Section~\ref{S5}).
In the broader setting of second-order operators acting on weighted spaces considered in \cite{Ebe} this implication
is not always valid.
The weights can introduce singular boundary behaviour (see \cite{Ebe}, Remark following Corollary~3.4).

Combination of the foregoing theorems gives the  conclusion that 
Markov uniqueness, $L_1$-uniqueness and the conservation property 
are all characterized by the capacity condition $\capp_\Omega(\partial\Omega)=0$.
This is of interest since the latter condition can be estimated in terms of the boundary
behaviour of the coefficients $c_{ij}$ and the geometric properties of $\partial\Omega$.
In Section~\ref{S4.2}  we derive estimates in terms of the order of degeneracy of
the coefficients and the Minkowski dimension of the boundary (see Proposition~\ref{psm4.1}).

The proofs of Theorems~\ref{tsm1.0},  \ref{tsm1.1}  and \ref{tsm1.2} will be given in Sections~\ref{S4.1}, \ref{S4.2} and \ref{S5}, respectively.
In Section~\ref{S6} we demonstrate that versions of the capacity estimates also give
sufficient conditions for $L_p$-uniqueness for all $p\in[1,2]$
and we establish that the semigroup $S^F_t$ is irreducible
if and only if $\Omega$ is connected.

\section{Elliptic regularity}\label{S3}

In this section we derive some basic regularity properties of  the operators $H_\Omega$ defined by~(\ref{ese1.0}).
 Since   $H_\Omega$ is symmetric its adjoint $H_\Omega^*$ is 
 an extension of its closure $\overline H_\Omega$ and 
the domain  $D(K_\Omega)$  of each self-adjoint extension $K_\Omega$   of $H_\Omega$
satisfies $D(\overline H_\Omega)\subseteq D(K_\Omega)\subseteq D(H_\Omega^*)$.
The principal observation is that  $D(\overline H_\Omega)$ and $ D(H_\Omega^*)$ only differ on the boundary $\partial\Omega$.
Hence the various possible extensions are distinguished by their boundary behaviour.

The comparison of $D(\overline H_\Omega)$ and $ D(H_\Omega^*)$ can be articulated in various ways but it is convenient for the sequel 
to express  it as a multiplier property.

\begin{thm}\label{tse3.1} Adopt the assumptions of Theorem~$\ref{tsm1.0}$.
Then
$C_c^\infty(\Omega)D(H_\Omega^*)\subseteq D(\overline H_\Omega)$.
\end{thm}
\proof\ 
The principal step in the proof consists of establishing that  $D(H_\Omega^*)\subseteq W^{1,2}_{\rm loc}(\Omega)$.
Once this is achieved the rest of the proof is given by the following argument.

Let $\Omega'$ be a bounded open subset of $\Omega$ which is strictly contained in $\Omega$, i.e.\
$\overline {\Omega'}\subset \Omega$.
(Strict containment will be denoted by $\Omega'\subset\subset\Omega$.)
If $\psi\in D(H_\Omega^*)$  and  $D(H_\Omega^*)\subseteq W^{1,2}_{\rm loc}(\Omega)$
then $\psi\in W^{1,2}(\Omega')$.
Set $\xi=H^*_\Omega\psi$ then $\xi\in L_2(\Omega')$ and
\[
\sum^d_{i,j=1}(c_{ij}\partial_j\eta,\partial_i\psi)=(\eta, \xi)
\]
for all $\eta\in C_c^\infty(\Omega')$, i.e.\  
$\psi$ is a weak solution
of the elliptic equation $H_{\Omega'}\psi=\xi$ on $\Omega'$.
Since $H_{\Omega'}$ is strongly elliptic on $L_2(\Omega')$  it  follows by elliptic regularity (see, for example, \cite{GT} Theorem~8.8) that  $\psi\in W^{2,2}(\Omega'')$ for all $\Omega''\subset\subset \Omega'$.
Thus $\psi\in W^{2,2}_{\rm loc}(\Omega')\subseteq  W^{2,2}_{\rm loc}(\Omega)$.
Therefore $D(H_\Omega^*)\subseteq  W^{2,2}_{\rm loc}(\Omega)$.
Hence if $\eta\in C_c^\infty(\Omega)$ then $\eta\,\psi\in W^{2,2}_0(\Omega)$.
But   $W^{2,2}_0(\Omega)\subseteq D(\overline H_\Omega)$, because the coefficients are bounded,
and the statement of the theorem is established.

It remains to prove that $D(H_\Omega^*)\subseteq W^{1,2}_{\rm loc}(\Omega)$.

First, fix $\eta\in C_c^\infty(\Omega)$ and set $K=\supp\eta$.
 Next  let   $\Omega'\subset\subset\Omega$ be a bounded open subset which contains $K$.
 Since $c_{ij}\in W^{1,\infty}(\Omega)$ and $C(x)>0$ for all  $x\in \Omega$  the restriction $H_{\Omega'}$ of $H_\Omega$ to $C_c^\infty(\Omega')$ is strongly elliptic.

Secondly,  $\varphi\in D(H_\Omega^*)$ if and only if  there is an $a>0$ such that 
 \[
 |(\varphi,H_\Omega\psi)|\leq a\,\|\psi\|_2
 \]
 for all $\psi\in C_c^\infty(\Omega)$.
 In particular if $\varphi\in D(H_\Omega^*)$ then these bounds are valid for all  $\psi\in C_c^\infty(\Omega')$.
 Thus the restriction   $\one_{\Omega'}\varphi$ of $\varphi$ to $\Omega'$  is in 
  $ D(H_{\Omega'}^*)$ and $H_{\Omega'}^*(\one_{\Omega'}\varphi)=
 \one_{\Omega'}(H_\Omega^*\varphi)$. 
 But 
$ \eta\,\varphi=\one_{\Omega'}(\eta\,\varphi)=\eta\,(\one_{\Omega'}\varphi)$.
In particular  if $\eta\,(\one_{\Omega'}\varphi)\in D(\overline H_{\Omega'})$ then 
 $ \eta\,\varphi\in D( \overline H_\Omega)$.
 Thus $\eta D(H^*_\Omega)\subseteq D(\overline H_\Omega)$ if and only if 
 $\eta D(H^*_{\Omega'})\subseteq D(\overline H_{\Omega'})$
  for all possible choices of $\eta$ and $\Omega'$.
 Therefore it suffices to prove $D(H_{\Omega'}^*)\subseteq W^{1,2}_{\rm loc}(\Omega')$
 for the strongly elliptic operator $H_{\Omega'}$ on $L_2(\Omega')$ for  all bounded open subsets $\Omega'\subset\subset\Omega$

Thirdly, we extend $H_{\Omega'}$ to a strongly elliptic operator $L$ on $L_2(\Ri^d)$ 
 with coefficients $\hat c_{ij}\in W^{1,\infty}(\Ri^d)$  such that $H_{\Omega'}=L|_{C_c^\infty(\Omega')}$.
This is achieved in two steps.
Since the $c_{ij}$ are continuous on $\Omega$, $C(x)\geq \mu I$ for all $x\in \Omega'$ and $C(x)>0$ for all $x\in \Omega$
one may choose an $\Omega''$ such that $\Omega'\subset\subset\Omega''\subset\subset\Omega$ and $C(x)\geq (\mu/2)I$ for all $x\in\Omega''$.
Then one may choose a $\chi\in C^\infty(\Ri^d)$ such that $0\leq \chi\leq 1$, $\chi(x)=1$ if $x\in\Omega'$ and $\chi(x)=0$ if $x$ is in the complement of $\Omega''$.
Then set $\widehat C=\chi \,C+(1-\chi)\,(\mu/2)I$.
It follows that $\widehat C\geq (\mu/2) I$.
Now let $L$ be the divergence form operator on $L_2(\Ri^d)$  with the matrix of coefficients $\widehat C=(\hat c_{ij})$.
It is strongly elliptic,  $\hat c_{ij}\in W^{1,\infty}(\Ri^d)$  and $L|_{C_c^\infty(\Omega')}=H_{\Omega'}$ by construction.
 Therefore the proof is completed by the following lemma.

\begin{lemma}\label{lse3.1} 
Let $H_{\Omega'}=L|_{C_c^\infty(\Omega')}$ where $L$ is a strongly elliptic operator, with  coefficients   $ \hat c_{ij}\in W^{1,\infty}(\Ri^d)$, acting on $L_2(\Ri^d)$.
Then  $D(H_{\Omega'}^*)\subseteq W^{1,2}_{\rm loc}(\Omega')$.
\end{lemma}
\proof\
The proof exploits some basic properties of strongly elliptic operators with Lipschitz continuous coefficients summarized in Proposition~\ref{pse1}
of the appendix. 
In particular $L$ is essentially self-adjoint on $C_c^\infty(\Ri^d)$ and its self-adjoint closure $\overline L$ has domain $D(\overline L)=W^{2,2}(\Ri^d)$.

Let $\cd({\Omega'})$ denote $C_c^\infty({\Omega'})$ equipped with the Frechet topology
and $\cd'({\Omega'})$ the dual space, i.e.\ the space of distributions on $\Omega'$.
If  $\psi\in L_2(\Omega')$ then $\varphi\in C_c^\infty (\Omega')\mapsto (\psi,H_{\Omega'}\varphi)$ is a continuous linear function over $\cd(\Omega')$.
Thus for each  $\psi\in L_2(\Omega')$ there is a  distribution $H_{\Omega'}(\psi)\in \cd'(\Omega')$ such that 
\[
(H_{\Omega'}(\psi),\varphi)=(\psi,H_{\Omega'}\varphi)
\]
for all $\varphi\in C_c^\infty(\Omega')$.
Similarly for each   $\psi\in L_2(\Ri^d)$ there is a  distribution $L(\psi)\in W^{-2,2}(\Ri^d)$, the dual of $W^{2,2}(\Ri^d)$, 
such that 
\[
(L(\psi),\varphi)=(\psi, L\varphi)
\]
for all $\varphi\in C_c^\infty(\Ri^d)$.
But by assumption
\[
(\psi,H_{\Omega'}\varphi)=(\psi, L\varphi)=(L(\psi),\varphi)
\]
for all $\psi\in L_2(\Omega')$ and all $\varphi\in C_c^\infty (\Omega')$.
Therefore 
\[
(H_{\Omega'}(\psi),\varphi)=(L(\psi),\varphi)
\]
for all $\varphi\in C_c^\infty (\Omega')$.
In particular $\psi\in D(H_{\Omega'}^*)$ if and only if $L(\psi)\in L_2(\Omega')$.

Next fix $\psi\in D(H_{\Omega'}^*)$.
Then  $L(\psi)\in L_2({\Omega'})$.
Moreover, if  $\eta\in C_c^\infty({\Omega'})$
then  $\eta\,\psi\in D(H_{\Omega'}^*)$ if and only if 
 $L(\eta\,\psi)\in L_2({\Omega'})$.
 But one has the distributional relation
\begin{equation}
L(\eta\,\psi)=\eta L(\psi)+L(\eta) \psi+\Psi_\eta
\label{ereg2.1}
\end{equation}
where 
\[
\Psi_\eta=-2\sum^d_{i,j=1}\hat c_{ij}\,(\partial_j\eta)\,(\partial_i\psi)
\;.
\]
Since $\psi, L(\psi)\in L_2({\Omega'})$
and $\eta, L(\eta) \in L_\infty(\Omega')$ it follows that the  first two terms on the right 
of (\ref{ereg2.1})  are  in $L_2({\Omega'})$.
It remains to demonstrate that $\Psi_\eta\in L_2({\Omega'})$.
But  there is an $a>0$ such that 
\[
|(\Psi_\eta,\varphi)|\leq a\,\|\psi\|_2\,\|\varphi\|_{1,2}
\]
for all $\varphi\in C_c^\infty(\Ri^d)$ where $\|\cdot\|_{1,2}$ denotes the $W^{1,2}$-norm. 
Therefore $\Psi_\eta\in W^{-1,2}(\Ri^d)$, the dual of $W^{1,2}(\Ri^d)$.
Hence  it follows from (\ref{ereg2.1})  that  $L(\eta\,\psi)\in W^{-1,2}(\Ri^d)$.
But 
\[
\eta\,\psi=(I+L)^{-1}\eta\,\psi+L(I+L)^{-1}\eta\,\psi=(I+L)^{-1}\eta\,\psi+(I+L)^{-1}L(\eta\,\psi)
\]
and $\eta\,\psi\in W^{1,2}(\Ri^d)$ by Proposition~\ref{pse1}.\ref{pse1-3} of the appendix applied to the strongly elliptic
operator $L$.
Since $\eta\in C_c^\infty({\Omega'})$ it follows that $\eta\,\psi\in W^{1,2}_0({\Omega'})$.

Finally let $K$ be an arbitrary compact subset of ${\Omega'}$.
If  $\eta_1\in C_c^\infty({\Omega'})$ with $\eta_1=1$ on $K$ it follows that 
 $\partial_i(\eta_1\psi)|_K=\partial_i\psi|_K$.
 Thus $\partial_i\psi\in L_2(K)$ for all $i\in\{1,\ldots,d\}$.
Therefore   $\psi\in W^{1,2}_{\rm loc}(\Omega)$.
\hfill$\Box$

\bigskip

One can also draw a conclusion about the domain of a general self-adjoint extension of~$H_\Omega$
and partially establish Statement~\ref{tsm1.0-3} of Theorem~\ref{tsm1.0}.

\begin{cor}\label{csme3.11} If  $K_\Omega$ is any  self-adjoint extension of $H_\Omega$ then 
$C^\infty_c(\Omega)D(K_\Omega)\subseteq D(\overline H_\Omega)$.
\end{cor}
\proof\
Since $D(\overline H_\Omega)\subseteq D(K_\Omega)\subseteq D(H^*_\Omega)$
one has $C^\infty_c(\Omega)D(K_\Omega)\subseteq C^\infty_c(\Omega)D(H^*_\Omega)\subseteq D(\overline H_\Omega)$ 
by Theorem~\ref{tse3.1}.
\hfill$\Box$

\bigskip

Note that $D(\overline H_\Omega)\subseteq D(H^*_\Omega)\subseteq W^{2,2}_{\rm loc}(\Omega)$.
Therefore $C^\infty_c(\Omega)D(K_\Omega)\subseteq 
D(\overline H_\Omega)\subseteq D(K_\Omega)\subseteq W^{2,2}_{\rm loc}(\Omega)$
for each self-adjoint extension $K_\Omega$.

\section{The minimal Markov extension}\label{S4.1}

In this section we prove Theorem~\ref{tsm1.0}.
The proof of the first parts  of the theorem broadly follows the reasoning used  in \cite{FOT} to prove 
the analogous result, Theorem~3.3.1,  for operators with $C^\infty$-coefficients.
The essential new ingredient is the elliptic regularity properties of Theorem~\ref{tse3.1}
and its corollaries.

\medskip

 \noindent{\bf Proof of Theorem~\ref{tsm1.0} }
 \noindent\ref{tsm1.0-0}.$\;$ The Markov property  of $l_\Omega$ follows by the calculations of Example~1.2.1 in \cite{FOT}.
 Moreover, the form is closed with respect to the graph norm by the arguments in Section~II.2.b of 
\cite{MR}.
The latter arguments depend crucially on the local strong ellipticity property.
Therefore the form $l_\Omega$ is a Dirichlet form. 
Finally   $D(l_\Omega)\cap C^\infty(\Omega)$ is a core of $l_\Omega$ by the proof of Lemma~3.3.3 in \cite{FOT}.

\medskip

 \noindent\ref{tsm1.0-1}.$\;$ 
The proof that  $L_\Omega$  is an extension of~$H_\Omega$ is identical to the proof of   Lemma~3.3.4 in \cite{FOT}
modulo a regularity argument.

Let $\varphi\in L_2(\Omega)$. 
Then $\psi=(I+L_\Omega)^{-1}\varphi\in D(L_\Omega)\subseteq D(l_\Omega)$.
Moreover,
\[
(\varphi,\eta)=(\psi,\eta)+l_\Omega(\psi,\eta)
=(\psi,\eta)+\int \sum^d_{i,j=1}c_{ij}\,(\partial_i\psi)(\partial_j\eta)
\]
for all $\eta\in C_c^\infty(\Omega)$.
Now fix an $\eta_1\in C_c^\infty(\Omega)$ such that $\eta_1=1$ on the support of $\eta$.
Then $\psi_1=\eta_1\psi\in D(l_\Omega)$ by a straightforward estimate and 
$\psi_1\in W^{1,2}_0(\Omega)$ by local strong ellipticity.
Therefore
\begin{eqnarray*}
(\varphi,\eta)
=(\psi_1,\eta)+\int\sum^d_{i,j=1} c_{ij}\,(\partial_i\psi_1)(\partial_j\eta)
=(\psi_1,(I+H_\Omega)\eta)=(\psi,(I+H_\Omega)\eta)
\end{eqnarray*}
 by partial integration.
Hence $\psi\in D(I+H^*_\Omega)$ and $(I+H^*_\Omega)\psi=\varphi$.
Thus $D(L_\Omega)\subseteq D(H^*_\Omega)$ and $ H_\Omega^*$ is an extension of $L_\Omega$.
So $D(H_\Omega)\subseteq D(L_\Omega)$ and $L_\Omega$ is an extension of $H_\Omega$.

\medskip

 \noindent\ref{tsm1.0-2}.$\;$ 
 This is the lengthiest part of the proof.
 We divide it into two steps.

\smallskip

\noindent{\bf Step 1}$\;$
First, we prove that $D(k_\Omega)\cap  D(H^*_\Omega)\subseteq D(l_\Omega)$
(see \cite{FOT}, proof of  Lemma~3.3.5).
Clearly it suffices to prove that 
\begin{equation}
k_\Omega(\varphi)\geq \int_\Omega dx\,  \sum^d_{i,j=1}c_{ij}(x)\,(\partial_i\varphi)(x)(\partial_j\varphi)(x)
\label{emse5.2}
\end{equation}
for all $\varphi\in D(k_\Omega)\cap D(H^*_\Omega)$.

Set $R_\lambda=(\lambda I+K_\Omega)^{-1}$ for all $\lambda>0$ 
and introduce the bounded forms
\begin{equation}
k^{(\lambda)}_\Omega(\varphi)=\lambda\,(\varphi,(I-\lambda R_\lambda)\varphi)
\label{emse5.3}
\end{equation}
for all  $\varphi\in L_2(\Omega)$.
The $k^{(\lambda)}_\Omega$
are Dirichlet forms and 
\begin{equation}
k_\Omega(\varphi)=\sup_{\lambda>0}k^{(\lambda)}_\Omega(\varphi)=\limsup_{\lambda\to\infty}k^{(\lambda)}_\Omega(\varphi)
\label{emse5.4}
\end{equation}
with $D(k_\Omega)$ the subspace of $L_2(\Omega)$ for which the supremum is finite (see, for example, \cite{FOT}, Lemma~1.3.4(ii)).

Next for $\eta\in C_c^\infty(\Omega)$ with $0\leq \eta\leq 1$ define the truncated form 
\[
k^{(\lambda)}_{\Omega,\eta}(\varphi)=k^{(\lambda)}_\Omega(\varphi,\eta\varphi)-2^{-1}k^{(\lambda)}_\Omega(\eta,\varphi^2)
\]
for all $\varphi\in L_2(\Omega)\cap L_\infty(\Omega)$.
It then follows from the Dirichlet form structure  that 
\begin{equation}
k^{(\lambda)}_\Omega(\varphi)\geq k^{(\lambda)}_{\Omega,\eta}(\varphi)
\label{emse5.6}
\end{equation}
for all $\varphi\in L_2(\Omega)\cap L_\infty(\Omega)$ (see \cite{BH}, Proposition~I.4.1.1).
Moreover,
\begin{equation}
k^{(\lambda)}_{\Omega,\eta}(\varphi)=
\lambda\,(\varphi, (I-\lambda R_\lambda)\eta\varphi)-
2^{-1}\lambda\,(\varphi(I-\lambda R_\lambda)\eta,\varphi)
\label{emse5.5}
\end{equation}
for all $\varphi\in L_2(\Omega)\cap L_\infty(\Omega)$ 
since $(I-\lambda R_\lambda)\eta\in L_\infty(\Omega)$ by the submarkovian property of $K_\Omega$.
Then, however, (\ref{emse5.5}) extends to all $\varphi\in L_2(\Omega)$ by continuity.
Combination of (\ref{emse5.3}), (\ref{emse5.4}), (\ref{emse5.6}) and (\ref{emse5.5}) immediately  gives
\begin{equation}
k_\Omega(\varphi)\geq \lambda\,(\varphi, (I-\lambda R_\lambda)\eta\varphi)-2^{-1}\lambda\,(\varphi(I-\lambda R_\lambda)\eta,\varphi)
\label{emse5.60}
\end{equation}
for all $\varphi\in D(k_\Omega)$.
Now we consider the limit $\lambda\to\infty$.

If $\varphi\in D(H_\Omega^*)$ then $\eta\varphi\in D(\overline H_\Omega)\subseteq D(K_\Omega)$
by Theorem~\ref{tse3.1}.
Therefore
\[
\lim_{\lambda\to\infty}\lambda\,(\varphi, (I-\lambda R_\lambda)\eta\varphi)=(\varphi, K_\Omega\eta\varphi)
=(\varphi, \overline H_\Omega\eta\varphi)
\;.
\]
Now let $S$ denote  the submarkovian semigroup generated by $K_\Omega$ on $L_2(\Omega)$ and $S^{(\infty)}$ the corresponding weak$^*$ semigroup on $L_\infty(\Omega)$.
Further let $K^{(\infty)}_\Omega$ denote  the generator of $S^{(\infty)}$ and $R^{(\infty)}_\lambda=(\lambda I+K^{(\infty)}_\Omega)^{-1}$ the resolvent.
Then 
$\eta\in D(H_\Omega)\cap L_\infty(\Omega)\subseteq D(K_\Omega)\cap L_\infty(\Omega)$ and $K_\Omega\eta=H_\Omega\eta\in L_\infty(\Omega)$.
Therefore $\eta\in D(K^{(\infty)}_\Omega)$ and $K^{(\infty)}_\Omega\eta=H_\Omega\eta$.
Consequently
\[
H_\Omega\eta=K^{(\infty)}_\Omega\eta={{\rm weak}^*\lim}_{\lambda\to\infty}\lambda(I-\lambda R^{(\infty)}_\lambda)\eta
\]
and one concludes that 
\[
\lim_{\lambda\to\infty}\lambda\,(\varphi(I-\lambda R_\lambda)\eta,\varphi)=
\lim_{\lambda\to\infty}\lambda((I-\lambda R^{(\infty)}_\lambda)\eta,\varphi^2)=(K^{(\infty)}_\Omega\eta,\varphi^2)=
(\varphi H_\Omega\eta,\varphi)
\;.
\]
Then it follows from taking the limit $\lambda\to\infty$ in (\ref{emse5.60}) that 
\[
k_\Omega(\varphi)\geq 
(\varphi, \overline H_\Omega\eta\varphi)-2^{-1}(\varphi H_\Omega\eta,\varphi)
\]
for all $\varphi\in D(k_\Omega)\cap D(H^*_\Omega)$.
Since $\eta\in C_c^\infty(\Omega)$ and $ D(H_\Omega^*)\subseteq W^{2,2}_{\rm loc}(\Omega)$, by the proof of Theorem~\ref{tse3.1},
it follows by direct calculation that
\[
k_\Omega(\varphi)\geq  \int_\Omega dx\, \sum^d_{i,j=1}\eta(x)\,c_{ij}(x)\,(\partial_i\varphi)(x)(\partial_j\varphi)(x)
\]
for all $\varphi\in D(k_\Omega)\cap D(H^*_\Omega)$ and $\eta\in C_c^\infty(\Omega)$ with $0\leq\eta\leq 1$.
But $k_\Omega(\varphi)$ is independent of~$\eta$.
Therefore taking the limit over a sequence of $\eta$ which converges monotonically upward to~$\one_\Omega$
one deduces  that (\ref{emse5.2}) is valid by  the Lebesgue dominated convergence theorem. 

\medskip
\noindent{\bf Step 2}$\;$
Next we argue that the inclusion $D(k_\Omega)\cap  D(H^*_\Omega)\subseteq D(l_\Omega)$
established by Step~1 implies $D(k_\Omega)\subseteq D(l_\Omega)$.

By definition $D(\overline h_\Omega)$ is a subspace of $ D(k_\Omega)$.
But the  orthogonal
complement of  $D(\overline h_\Omega)$  with respect to the graph norm 
$\|\cdot\|_{ D(k_\Omega)}$
is  $ D(k_\Omega)\cap \cn$
where 
\[
\cn=\{\varphi\in D(H_\Omega^*): (I+H^*_\Omega)\varphi=0\}
\]
(see \cite{FOT}, Lemma~3.3.2(ii)).
Therefore each $\varphi\in D(k_\Omega)$ has a 
unique decomposition $\varphi=\varphi_1+\varphi_2$ with $\varphi_1\in D(\overline h_\Omega)$ and $\varphi_2\in D(k_\Omega)\cap \cn$ such that 
\[
\|\varphi\|_{ D(k_\Omega)}^2=
\|\varphi_1\|_{ D(\overline h_\Omega)}^2+\|\varphi_2\|_{ D(k_\Omega)}^2
\;.
\]
But  $\varphi_2\in D(k_\Omega)\cap D(H^*_\Omega)$.
So $\varphi_2\in D(l_\Omega)$ and $k_\Omega(\varphi_2)\geq l_\Omega(\varphi_2)$  by Step~2.
Therefore
\begin{eqnarray*}
\|\varphi\|_{ D(k_\Omega)}^2\geq
\|\varphi_1\|_{ D(\overline h_\Omega)}^2+\|\varphi_2\|_{ D(l_\Omega)}^2
=\|\varphi_1\|_{ D(l_\Omega)}^2+\|\varphi_2\|_{ D(l_\Omega)}^2
=\|\varphi\|_{ D(l_\Omega)}^2
\;.
\end{eqnarray*}
The last equality follows because $l_\Omega(\varphi_1,\varphi_2)+(\varphi_1,\varphi_2)=0$.
It follows 
immediately that $D(k_\Omega)\subseteq D(l_\Omega)$.
This completes the proof of Statement~\ref{tsm1.0-2} of Theorem~\ref{tsm1.0}.

\medskip

 \noindent\ref{tsm1.0-3}.$\;$ 
The inclusion  
$C_c^\infty(\Omega)D(K_\Omega)
\subseteq D(\overline H_\Omega)$,
was established in  Corollary~\ref{csme3.11}.
But if $\eta\in C_c^\infty(\Omega)$ and $\varphi\in D(K_\Omega)$ then $\eta\varphi\in D(\overline H_\Omega)\subseteq D(K_\Omega)\subseteq D(k_\Omega)\subseteq D(l_\Omega)$ and 
\begin{eqnarray*}
\overline h_\Omega(\eta\varphi)=l_\Omega(\eta\varphi)\leq 2\,\|\eta\|_\infty^2\,l_\Omega(\varphi)+2\,\|\Gamma(\eta)\|_\infty\,\|\varphi\|_2^2
\leq 2\,(\|\Gamma(\eta)\|_\infty+\|\eta\|_\infty^2)\,\|\varphi\|_{D(k_\Omega)}^2
\end{eqnarray*}
where $\Gamma(\eta)=\sum^d_{i,j=1}c_{ij}\,(\partial_i\eta)\,(\partial_j\eta)$,
i.e.\ $\Gamma$ is the {\it carr\'e du champ} as defined in  \cite{BH}, Section~I.8.
Since $D(K_\Omega)$ is a core of $k_\Omega$ with respect to the $D(k_\Omega)$-graph norm 
it follows that $C_c^\infty(\Omega)D(k_\Omega)
\subseteq D(\overline h_\Omega)$ by continuity.

\medskip

 \noindent\ref{tsm1.0-5}.$\;$ First let $D(k_\Omega)_c$ denote the subspace of functions with compact support in $D(k_\Omega)$.
If $\varphi\in D(k_\Omega)_c$ then  by regularization one can construct a sequence $\varphi_n\in C_c^\infty(\Omega)$ such that 
$\|\varphi_n-\varphi\|_{D(k_\Omega)}\to0$ as $n\to\infty$.
Since $k_\Omega$ is an extension of $h_\Omega$ it follows that $\varphi\in D(\overline h_\Omega)$.
Therefore $D(k_\Omega)_c\subseteq D(\overline h_\Omega)$.
Now the first part of Statement~\ref{tsm1.0-5} follows from Proposition~2.1 of \cite{ER30}.
But then $D(\overline h_\Omega)$ is an ideal (see \cite{Ouh5}, Definition~2.19) of $D(k_\Omega)$ by Corollary~2.22  of  \cite{Ouh5}.
In particular it is an order ideal.

For the proof of the second part of Statement~\ref{tsm1.0-5} we   again appeal to Corollary~2.22  of  \cite{Ouh5}.
First if  $e^{-tL_\Omega}$ dominates $ e^{-tK_\Omega}$ then it follows from this corollary that $D(k_\Omega)$ is an ideal of $D(l_\Omega)$.
Secondly, for the converse statement,  it  suffices to prove that $D(k_\Omega)$ is an ideal of $D(l_\Omega)$.
Then the domination property follows from another application of Corollary~2.22 of \cite{Ouh5}.
Thus if  $\psi\in D(k_\Omega)$, $\varphi\in D(l_\Omega)$  and $|\varphi|\leq |\psi|$ then 
one  must deduce  that $\varphi\sgn\psi\in D(k_\Omega)$.
But $\varphi, \psi\in D(l_\Omega)$. 
Therefore  $\varphi\sgn\psi\in D(l_\Omega)$ by Proposition~2.20 of \cite{Ouh5}.
(See the remark following this proposition.)
Moreover,  $|\psi|\in D(k_\Omega)$ and $(\varphi\sgn\psi)_+\in D(l_\Omega)$      because $k_\Omega$ and $l_\Omega$  are Dirichlet forms.
Since
\[
0\leq (\varphi\sgn\psi)_+\leq|\varphi|\leq|\psi|
\]
and since $D(k_\Omega)$ is an order ideal of $D(l_\Omega)$ it follows  that $ (\varphi\sgn\psi)_+\in D(k_\Omega)$.
Applying the same argument to $-\varphi$ one deduces that $ (\varphi\sgn\psi)_-\in D(k_\Omega)$.
Therefore $\varphi\sgn\psi\in D(k_\Omega)$ and $D(k_\Omega)$ is an ideal of $D(l_\Omega)$.
\hfill$\Box$

\bigskip

We note in passing that the existence of   $l_\Omega$ gives a criterion for uniqueness of the submarkovian extension of $H_\Omega$
similar to the standard criterion for essential self-adjointness.

\begin{prop}\label{psm4.0}
The following conditions are equivalent.
\begin{tabel}
\item\label{psm4.0-1}
$H_\Omega$ has a unique submarkovian extension.
\item\label{psm4.0-2}
$\;\ker(I+H^*_\Omega)\cap D(l_\Omega)=\{0\}$
\end{tabel}
\end{prop}
\proof\
 Condition~\ref{psm4.0-1} is equivalent to $l_\Omega=\overline h_\Omega$
 by Theorem~\ref{tsm1.0}, i.e.\ equivalent to $D(l_\Omega)=D(\overline h_\Omega)$.
 But $D(l_\Omega)=D(\overline h_\Omega)\oplus \ch_\Omega$ with $\ch_\Omega=\ker(I+H^*_\Omega)\cap D(l_\Omega)$
by  Lemma~3.3.2(ii) in \cite{FOT}.
Therefore the equivalence of the conditions of the proposition  is immediate.
\hfill$\Box$

\bigskip

Statement~\ref{tsm1.0-2} of Theorem~\ref{tsm1.0} establishes that $H_\Omega$ is Markov unique if and only if 
$l_\Omega=\overline h_\Omega$ or, equivalently, $D(l_\Omega)=D(\overline h_\Omega)$.
This is the criterion  used extensively in the analysis of Markov uniqueness (see \cite{Ebe}, Chapter~3).
It  will also be used to prove Theorem~\ref{tsm1.1}. 

Statement~\ref{tsm1.0-3} of the theorem establishes  that each submarkovian extension of $H_\Omega$ is a Silverstein extension
(see \cite{Take} or  \cite{Ebe}, Definition~1.4).
Therefore Markov uniqueness of $H_\Omega$ and Silverstein uniqueness are equivalent.

Statement~\ref{tsm1.0-5} gives an alternative approach to establishing Markov uniqueness of $H_\Omega$ if the 
submarkovian semigroup generated by the Friedrichs extension is conservative.
This will be discussed in Section~\ref{S5}.

\section{Markov Uniqueness}\label{S4.2}

In this section we prove Theorem~\ref{tsm1.1}.
Throughout the section we assume that  the coefficients $c_{ij}$  are real, symmetric, Lipschitz continuous and   $C(x)=(c_{ij}(x))>0$
for all $x\in\Omega$.

\medskip

\noindent{\bf Proof of Theorem~\ref{tsm1.1} }
It follows from Theorem~\ref{tsm1.0}  that $H_\Omega$ has a unique submarkovian extension
if and only if $l_\Omega=\overline h_\Omega$, i.e.\ if and only if $C_c^\infty(\Omega)$ is a core of $l_\Omega$. 
Therefore Theorem~\ref{tsm1.1} is a direct corollary of the following proposition
which is a variation of Proposition~3.2 in \cite{RSi}.
  
  \begin{prop}\label{psm4.10}
  Under the assumptions of Theorem~$\ref{tsm1.1}$, the following conditions are equivalent:
  \begin{tabel}
 \item\label{psm4.10-1}
  $\capp_{\Omega}(\partial\Omega)=0$.
  \item\label{psm4.10-3}
 $C_c^\infty(\Omega)$ is a core of $l_\Omega$.
    \end{tabel}
  \end{prop}
  \proof 
\noindent  \ref{psm4.10-1}$\Rightarrow$\ref{psm4.10-3}$\;$
First, since $l_\Omega$ is a Dirichlet form $D(l_\Omega)\cap L_\infty(\Omega)$ is a core of $l_\Omega$.
  Therefore  it suffices  that each $\varphi\in D(l_\Omega)\cap L_\infty(\Omega)$ can be approximated by a sequence $\varphi_n\in C_c^\infty(\Omega)$
  with respect to the graph norm $\|\cdot\|_{D(l_\Omega)}$.
  Now fix $\varphi\in D(l_\Omega)\cap L_\infty(\Omega)$.

Secondly, let $\rho_n\in C_c^\infty(\Ri^d)$ be a sequence of functions with $0\leq\rho_n\leq1$, $\|\nabla\rho_n\|_\infty\leq {n^{-1}}$  and such that $\rho_n\to\one$ pointwise as
$n\to\infty$.
Then $\one-\rho_n\in W^{1,\infty}(\Ri^d)$.
But $W^{1,\infty}(\Ri^d)D(l_\Omega)\subseteq D(l_\Omega)$.
Therefore $(\one-\rho_n)\varphi\in D(l_\Omega)\cap L_\infty(\Omega)$.
It then  follows from Leibniz' rule and the Cauchy--Schwarz inequality that 
\begin{eqnarray*}
\|\varphi-\rho_n\varphi\|_{D(l_\Omega)}^2&\leq &2\,\int_\Omega \varphi^2\,\Gamma(\rho_n)+2\int_\Omega(\one-\rho_n)^2\,\Gamma(\varphi)+\|(\one-\rho_n)\varphi\|_2^2\\
&\leq&2\,n^{-2}\|C\|\,\|\varphi\|_2^2+\int_\Omega(\one-\rho_n)^2\,(2\,\Gamma(\varphi)+\varphi^2)
\end{eqnarray*}
where $\|C\|$ is the supremum over the matrix norms $\|C(x)\|$.
Clearly the first term on the right hand side tends to zero as $n\to\infty$.
Moreover,  $0\leq (\one-\rho_n)^2\leq1$, $(\one-\rho_n)^2\to0$  pointwise as $n\to\infty$ and $2\,\Gamma(\varphi)+\varphi^2\in L_1(\Omega)$.
Therefore the second term on the right hand side also tends to zero as $n\to\infty$ by the Lebesgue dominated convergence theorem.
Thus $\varphi$ is approximated by the sequence $\rho_n\varphi$ in the graph norm.

Thirdly, since $\capp_{\Omega}(\partial\Omega)=0$ one may choose    $\chi_n\in D(l_\Omega)\cap L_\infty(\Omega)$ and open subsets 
$U_n\supset \partial\Omega$ such that $0\leq \chi_n\leq 1$, $l_\Omega(\chi_n)+\|\chi_n\|_2^2\leq n^{-1}$ and $\chi_n= 1$ on $U_n\cap\,\Omega$.
Now set $\varphi_{n}=(1-\chi_n)\rho_n\varphi$.
Then 
\[
\|\varphi-\varphi_{n}\|_{D(l_\Omega)}^2\leq 2\,\|\varphi-\rho_n\varphi\|_{D(l_\Omega)}^2+2\, \|\chi_n\rho_n \varphi\|_{D(l_\Omega)}^2
\]
and the first term on the right hand side converges to zero as $n\to\infty$ by the previous discussion.
Moreover,
\begin{equation}
\|\chi_n\rho_n \varphi\|_{D(l_\Omega)}^2= l_\Omega(\chi_n\rho_n\varphi)+\|\chi_n\rho_n\varphi\|_2^2
  \label{ebdn1}
\end{equation}
and the second term on the right tends to zero because $\|\chi_n\rho_n\varphi\|_2\leq \|\chi_n\|_2\|\varphi\|_\infty$.
But the first  term can be estimated by
\begin{eqnarray*}
l_\Omega(\chi_n\rho_n\varphi)
&\leq& 2\int_\Omega\rho_n^2\varphi^2\,\Gamma(\chi_n)+ 2\int_\Omega\chi_n^2\,\Gamma(\rho_n\varphi)\\
&\leq&2\int_\Omega\varphi^2\,\Gamma(\chi_n)+4\int_\Omega\chi_n^2\,\varphi^2\,\Gamma(\rho_n)+4\int_\Omega\chi_n^2\,\rho_n^2\,\Gamma(\varphi)\\
&\leq&2\,\|\varphi\|_\infty^2\,l_\Omega(\chi_n)+4\,\|C\|\,\|\nabla\rho_n\|_\infty^2\|\varphi\|_2^2+4\int_\Omega\chi_n^2\,\Gamma(\varphi)
\;.
\end{eqnarray*}
Since $l_\Omega(\chi_n)\to0$ and $\|\nabla\rho_n\|_\infty\to0$ as $n\to \infty$ the first two terms on the right hand side tend to zero.
But  if $A_m=\{x\in\Omega: \Gamma(\varphi)>m\}$ one has the equicontinuity estimate
\[
\int_\Omega\chi_n^2\,\Gamma(\varphi)\leq m\,\|\chi_n\|_2^2 +\int_{A_m}\Gamma(\varphi)
\]
because $0\leq \chi_n\leq 1$.
Since $\|\chi_n\|_2\to0$  and $\Gamma(\varphi)\in L_1(\Omega)$ the integral also tends to zero as $n\to\infty$.
Thus both  terms on the right hand side of (\ref{ebdn1}) tend to zero as $n\to\infty$ and 
one now concludes that $\varphi$ is approximated by the sequence $\varphi_{n}$ in the graph norm.

Finally $\supp\varphi_n$ is contained in the set $\Omega_n=((\supp\rho_n)\cap\,\Omega)\cap (\Omega\backslash(U_n\cap\,\Omega))$.
Hence   $\Omega_n$ is a bounded subset  which is strictly contained in $\Omega$,  i.e.\ $\Omega_n\subset\subset\Omega$. 
Then  since $C(x)>0$ for all $x\in\Omega$ one has an estimate
$l_\Omega(\varphi_n)\geq\mu_n\,\|\nabla\varphi_n\|_2^2$ with $\mu_n>0$. 
Therefore $\varphi_n\in W^{1,2}_0(\Omega_n)$ and it follows that it can be approximated in the $W^{1,2}(\Omega_n)$-norm by a sequence of $C_c^\infty$-functions.
Then because $l_\Omega(\psi)\leq \|C\|\,\|\nabla\psi\|_2^2$ for all $\psi\in D(l_\Omega)$ it follows that $\varphi_n$, and hence $\varphi$, can be approximated
by a sequence of $C_c^\infty$-functions in the graph norm $\|\cdot\|_{D(l_\Omega)}$.

\smallskip

\noindent \ref{psm4.10-3}$\Rightarrow$\ref{psm4.10-1} $\;$
Let $\psi\in D(l_\Omega)\cap C^\infty(\Omega)$ with $\psi=1$ on $U\cap\Omega$ where $U$ is an 
open subset containing $\partial\Omega$.
One may assume $0\leq \psi\leq 1$.
Then by Condition~\ref{psm4.10-3}
there is a sequence $\psi_n\in C_c^\infty(\Omega)$ such that $\|\psi_n-\psi\|_{D(l_\Omega)}\to0$.
In particular $\psi\in W^{1,2}(\Omega)$.
Since $\psi_n$ has compact support in $\Omega$ it also follows that there is an open subset containing $\partial\Omega$
such that $\psi_n=0$ on $U_n\cap\Omega$.
Therefore $\psi-\psi_n=1$ on $(U\cap U_n)\cap \Omega$ and one must have $\capp_{\Omega}(\partial\Omega)=0$.
\hfill$\Box$

\begin{remarkn}\label{rsm4.1}
Although $\capp_\Omega$ is defined in terms of the space $D(l_\Omega)$ it suffices for the discussion of Markov uniqueness
to consider the restriction to $\Omega$ of functions in $W^{1,2}(\Ri^d)$.
If, for example, one defines
\begin{eqnarray*}
\widetilde \capp_\Omega(A)=\inf\Big\{\|\psi\|_{W^{1,2}(\Omega)}^2&&: \psi\in W^{1,2}(\Ri^d) \mbox{ and  there exists   an open set  }\nonumber\\[-5pt]
&& U\subset \Ri^d
\mbox{ such that } U\supseteq A
\mbox{ and } \psi=1 \mbox{ a.\ e.\ on } U\cap\Omega\Big\}
\end{eqnarray*}
then $\widetilde \capp_\Omega(A)\geq \capp_\Omega(A)$ for all measurable subsets $A\subseteq \overline\Omega$.
Therefore  $\widetilde \capp_\Omega(\partial\Omega)=0$ implies that  $ \capp_\Omega(\partial\Omega)=0$.
Conversely, if  $ \capp_\Omega(\partial\Omega)=0$ then it follows from Proposition~\ref{psm4.10}
that $C^\infty_c(\Omega)$ is a core of $l_\Omega$.
Therefore the argument used to establish that  \ref{psm4.10-3}$\Rightarrow$\ref{psm4.10-1}  in the proof of Proposition~\ref{psm4.10}
also establishes that $\widetilde \capp_\Omega(\partial\Omega)=0$.
Thus  $ \capp_\Omega(\partial\Omega)=0$ if and only if $\widetilde \capp_\Omega(\partial\Omega)=0$.
\end{remarkn}

The condition $\capp_\Omega(\partial\Omega)=0$ is of interest as a criterion for Markov uniqueness since the capacity can be estimated
by elementary means.
The estimates depend on two gross features, the order of degeneracy of the coefficients at the boundary
and the dimension of the boundary.

Let $A$ be a measurable subset of $\overline\Omega$ with $|A|=0$
and $B$ a bounded measurable subset of $A$.
Then for each $\delta>0$ define the $\delta$-parallel body $B_\delta$ of $B$ by
\[
B_\delta=\{x\in \Ri^d:\inf_{y\in B}|x-y|<\delta\}
\;.
\]
There are a variety of ways of assigning a dimension to $A$ or $B$ (see, for example, \cite{Fal}, Chapters~2 and 3).
The Minkowski dimension $d(B)$ of $B$  is the smallest positive real value for which there is a $b>0$
such that $|B_\delta|\leq b\, \delta^{d-d(B)}$ for all $\delta\in\langle0,1]$.
In general $d(B)\in[0, d\rangle$ with $d(B)=0$ if $B$ is a finite set, $d(B)=1$ if $B$ is a line segment etc.
The dimension of $A$ is defined by $d(A)=\sup\{d(B): B\subseteq A\}$.

\begin{prop}\label{psm4.1}
 Let $A$ be a  measurable subset of   $\overline\Omega$.
Assume there are $a>0$ and $\gamma\geq0$ such that $0< C(x)\leq a\,(d_A(x)\wedge 1)^\gamma$ for all $x\in\Omega$
where $d_A$ is the Euclidean distance to $A$.

If $\gamma\geq 2-(d-d(A))$ then $\capp_\Omega(A)=0$.
In particular $\capp_\Omega(A)=0$ for all $A$ with $d(A)\leq d-2$.
\end{prop}
\proof\
 Let $B$ be a  bounded measurable subset of   $A$.
 If $\delta$ is sufficiently small then $d_A(x)=d_{B_\delta}(x)$ for all $x\in B_\delta$.
 For convenience choose units so that $\delta=1$.
 Then introduce the functions $x>0\mapsto \chi_n(x)\in[0,1]$ by $\chi_n(x)=1$ if $x\in\langle0,n^{-1}]$, $\chi_n(x)=-\log x/\log n$ if $x\in\langle n^{-1},1]$ and $\chi_n(x)=0$ if $x>1$.
Set $\eta_n=\chi_n\circ d_B$.
It follows that $\eta_n(x)=0$ if $x\not\in B_1$.
Now  $\|\eta_n\|_2\to0$ as $n\to\infty$.
 Moreover,
 \[
l_\Omega(\eta_n)\leq  a\int_{B_1} dx\, d_B(x)^{\gamma}|\nabla\eta_n(x)|^2\leq a\, (\log n)^{-2}\int dx\,\one_{\{x:n^{-1}\leq d_B(x)\leq1\}}\,d_B(x)^{-(2-\gamma)}
 \;.
 \]
Using the identity $d_B(x)^{-(2-\gamma)}=1+(2-\gamma)^{-1}\int^1_{d_B(x)}dt\,t^{-(3-\gamma)}$ 
 and changing the order of integration one immediately deduces that  
 \[
l_\Omega(\eta_n)\leq a'(\log n)^{-2}\Big(1+(2-\gamma)^{-1}\int_{n^{-1}}^1dt\,t^{-3+\gamma+d-d(B)}(t^{-(d-d(B))}|B_t|)\Big)
\;.
\]
Thus if $\gamma\geq  2-(d-d(B))$  then
$l_\Omega(\eta_n)\leq a''(\log n)^{-1}\to0$ as $n\to\infty$.
It follows that  $\capp_\Omega(B)=0$.
Then, however,  it follows from the general additivity properties of the capacity that $\capp_\Omega(A)=0$ for $\gamma\geq  2-(d-d(A))$ .
\hfill$\Box$

\bigskip

The estimates of Proposition~\ref{psm4.1} have two simple implications.

\begin{cor}\label{ced1}
Assume  $d(\partial\Omega)=d-1$.
If the coefficients $c_{ij}\in W^{1,\infty}_0(\Omega)$ are real symmetric
and  $C(x)>0$ for all $x\in \Omega$ then $H_\Omega$ is Markov unique.
\end{cor}
\proof\
Since the  coefficients $c_{ij}$ are in $ W^{1,\infty}_0(\Omega)$ they  extend by  continuity  
to $\overline\Omega$ and the extended coefficients are zero on the boundary.
But then by Lipschitz continuity $|c_{ij}(x)|\leq a\,(d_{\partial\Omega}(x)\wedge 1)$ for all $x\in\Omega$.
Therefore $\capp_\Omega(\partial\Omega)=0$ by Proposition~\ref{psm4.1} applied with $\gamma=1$.
Hence $H_\Omega$ is Markov unique by Theorem~\ref{tsm1.1}.
\hfill$\Box$

\bigskip

\begin{cor}\label{ced2}
If the coefficients $c_{ij}\in W^{2,\infty}_0(\Omega)$ are real symmetric
and  $C(x)>0$ for all $x\in \Omega$ then $H_\Omega$ is Markov unique.
\end{cor}
\proof\
Since the coefficients $c_{ij}$ are in $ W^{2,\infty}_0(\Omega)$ they again extend to $\overline\Omega$, the extensions are
 zero on the boundary and  one now  has bounds 
 $|c_{ij}(x)|\leq a\,(d_{\partial\Omega}(x)\wedge 1)^2$ for all $x\in\Omega$.
 Then $\capp_\Omega(\partial\Omega)=0$, by Proposition~\ref{psm4.1} applied with $\gamma=2$,
and $H_\Omega$ is Markov unique by Theorem~\ref{tsm1.1}.
\hfill$\Box$

\bigskip

Note that the second result is universal in the sense that it does not depend on the geometry of $\Omega$.
In particular it does not depend on the dimension of $\partial\Omega$.
Moreover, it suffices that the coefficients $c_{ij}\in W^{1,\infty}(\Omega)\cap W^{2,\infty}_0(U\cap\Omega)$ 
for some open set $U\supset\partial\Omega$.
In fact if  $c_{ij}\in W^{2,\infty}_0(\Omega)$ then the weaker ellipticity condition $C(x)\geq0$ for $x\in\Omega$
suffices to deduce that $H_\Omega$ is $L_2$-unique (see \cite{Rob7}, Section~6, or \cite{ER27}, Proposition~2.3).
In this latter case the coefficients can be extended to $\Ri^d$ by setting $c_{ij}(x)=0$ if $x\in\Omega^{\rm c}$ and 
then the operator is essentially self-adjoint  on $C_c^\infty(\Ri^d)$ and the self-adjoint extension leaves $L_2(\Omega)$
invariant.

Finally we emphasize  that the condition $\capp_\Omega(\partial\Omega)=0$ does not necessarily imply that
the coefficients  $c_{ij}(x)\to0$ as $x\to\partial\Omega$.
In fact Proposition~\ref{psm4.1} establishes that if $A\subset\partial \Omega$ and $d(A)\leq d-2$ then $\capp_\Omega(A)=0$ independently of the 
boundary behaviour of the coefficients.
Nevertheless if the boundary $\partial\Omega$ is Lipschitz continuous, and consequently $d(\partial\Omega)=d-1$,
one can argue that $\capp_\Omega(\partial\Omega)=0$ if and only if  $c_{ij}(x)\to0$ as $x\to\partial\Omega$.

\section{Conservation criteria}\label{S5}

In this section we prove  Theorem~\ref{tsm1.2}.
This theorem is to a large extent known and we concentrate on 
the new feature, Markov uniqueness implies semigroup conservation.
An integral part  in this proof is played by an approximation criterion
for conservation which is also useful for the discussion of $L_p$-uniqueness (see Section~\ref{S6}).

\begin{lemma}\label{lsm5.1}
Assume  there exists  a sequence $\eta_n\in C_c^\infty(\Omega)$ 
with $0\leq \eta_n\leq \one_\Omega$ such that $\|(\eta_n-\one_\Omega)\psi\|_2\to0$ for all $\psi\in L_2(\Omega)$  
and $h_\Omega(\eta_n)\to0$ as $n\to\infty$.
Then $S^F_t$ is conservative.
\end{lemma}
\proof\
First it follows that $((\eta_n-\one_\Omega),\psi)\to0$ as $n\to\infty$ for all $\psi\in L_1(\Omega)\cap L_2(\Omega)$.
Fix  $\varphi$ in the $L_1$-dense  set $D(H^F_\Omega)\cap L_1(\Omega)$.
Then $S^F_t\varphi\in L_1(\Omega)\cap L_2(\Omega)$ and
\begin{eqnarray*}
|(\one_\Omega, S^F_t\varphi)-(\one_\Omega,\varphi)|&=&\lim_{n\to\infty}|(\eta_n, S^F_t\varphi)-(\eta_n,\varphi)|\\
&=&\lim_{n\to\infty}\Big|\int^t_0ds\,(\eta_n,S^F_sH_\Omega^F\varphi)\Big|
\leq \lim_{n\to\infty}t\,h_\Omega(\eta_n)^{1/2}\,\overline h_\Omega(\varphi)^{1/2}=0
\;.
\end{eqnarray*}
Therefore $S^F_t$ is conservative on $L_\infty(\Omega)$.
\hfill$\Box$

\bigskip

Now we turn to the proof of the theorem

\smallskip

\noindent{\bf Proof of Theorem~\ref{tsm1.2}}$\;$
\noindent\ref{tsm1.2-3}$\Rightarrow$\ref{tsm1.2-0}$\;$
The proof is in five steps.

\smallskip

\noindent{\bf Step~1}$\;$ The first step consists of proving  the implication for bounded  $\Omega$  by 
constructing  a sequence of $\eta_n$ of the type occurring  in Lemma~\ref{lsm5.1}.
 
 Assume $\Omega$ is bounded.
 It follows from the  Markov uniqueness    and Theorem~\ref{tsm1.1} that $\capp_\Omega(\partial\Omega)=0$
 and $l_\Omega=\overline h_\Omega$.
Therefore  there exist
  a  decreasing sequence of  open subsets $U_n$ of $\Ri^d$ with $\partial\Omega\subset \,U_n$ and a sequence $\chi_n\in D(l_\Omega)$ with $0\leq \chi_n\leq 1$ and $\chi_n=1$  on $U_n\cap\Omega$
such that $\|\chi_n\|_2\to0$ and $l_\Omega(\chi_n)\to 0$  as $n\to\infty$.
Since $\Omega$ is bounded it follows that $\one_\Omega\in D(l_\Omega)$.
Therefore $\eta_n=(\one_\Omega-\chi_n)\in D(l_\Omega)$.
But then 
\[
\|(\eta_n-\one_\Omega)\psi\|_2=\|\chi_n\psi\|_2\leq \|\chi_n\|_2\,\|\psi\|_\infty\to0
\;.
\]
for all  $\psi$ in the $L_2$-dense subset $ L_2(\Omega)\cap L_\infty(\Omega)$.
Thus the first convergence property of the $\eta_n$ is satisfied.
Then, however, $l_\Omega(\eta_n)=l_\Omega(\chi_n)$ and the second condition is also satisfied.
Finally $\supp\eta_n\subset\subset\Omega$ for each $n$.
Hence by regularization  one may construct a second sequence of $C^\infty_c(\Omega)$-functions $\eta_n\in D(l_\Omega)$ with similar boundedness and convergence
properties.

Therefore it follows from  Lemma~\ref{lsm5.1} that  the semigroup $S^F_t$ is conservative.

\smallskip

\noindent{\bf Step~2}$\;$
The second step consists of proving the theorem for unbounded $\Omega$  but for a  family of cutoff operators.

Fix  $\rho\in C^\infty_c(\Ri^d)$ with $0\leq\rho\leq 1$, $\rho(x)=1$ if $|x|\leq 1$ and $\rho(x)=0$ if $|x|\geq 2$.
Then  introduce the sequence $\rho_n$ by $\rho_n(x)=\rho(n^{-1}x)$.
Thus $\rho_n(x)=1$ if $|x|\leq n$ and $\rho(x)=0$ if $|x|\geq 2n$.
Set $B_n=\{x\in\Ri^d:|x|<2n\}$ and $\Omega_n=\Omega\cap B_n$.
Note that $\Omega_n$ is bounded.

Now define a family of truncations $h_{\Omega,n}$ of $h_\Omega$ by $D(h_{\Omega, n})=C_c^\infty(\Omega_n)$
and
\[
h_{\Omega,n}(\varphi)=h_\Omega(\varphi, \rho_n\varphi)-2^{-1}h_\Omega(\rho_n,\varphi^2)
\]
for all $\varphi\in C_c^\infty(\Omega_n)$.
The truncation $h_{\Omega,n}$ is the Markovian form corresponding to the symmetric operator with $H_{\Omega,n}$ coefficients $\rho_nc_{ij}$ acting on $L_2(\Omega_n)$.
Let $l_{\Omega, n}$  denote the extended form corresponding to $H_{\Omega,n}$.
The form $h_{\Omega,n}$ is  automatically closable, the closure $\overline h_{\Omega,n}$ is a Dirichlet form and the corresponding  self-adjoint operator $H_{\Omega,n}^F$ is the  Friedrichs extension of $H_{\Omega,n}$.
The form $l_{\Omega, n}$ is a Dirichlet form which in principle differs from $\overline h_{\Omega,n}$.
But we next argue that $H_{\Omega,n}$ is Markov unique.
Hence $l_{\Omega,n}=\overline h_{\Omega,n}$.  

Let $\capp_{\Omega,n}(A)$ denote the capacity of the measurable subset $A$ of $\overline\Omega_n$ measured with respect to $l_{\Omega,n}$.
Since $\Omega_n=\Omega\cap B_n$ it follows that  $\partial\Omega_n=(\partial\Omega\cap \overline B_n)\cup(\partial B_n\cap\overline\Omega)$.
Hence
\[
\capp_{\Omega,n}(\partial\Omega_n)=\capp_{\Omega,n}(\partial\Omega\cap \overline B_n)+\capp_{\Omega,n}(\partial B_n\cap \overline\Omega)
\;.
\]
But $l_{\Omega, n}\leq l_\Omega$ and $\capp_\Omega(\partial\Omega)=0$ by  Markov uniqueness of $H_\Omega$ .
Therefore
\[
\capp_{\Omega,n}(\partial\Omega\cap \overline B_n)\leq \capp_{\Omega}(\partial\Omega\cap \overline B_n)\leq \capp_{\Omega}(\partial\Omega)=0
\;.
\]
Moreover, $\capp_{\Omega,n}(\partial B_n\cap \overline\Omega)=0$ because the $C^\infty$-cutoff function $\rho_n$ and all its derivatives
are zero on the boundary $\partial B_n$.
Thus $\capp_{\Omega,n}(\partial\Omega_n)=0$ and $H_{\Omega,n}$ is Markov unique by Theorem~\ref{tsm1.1}.
Hence the semigroup generated by the Friedrichs extension $H_{\Omega,n}^F$ of the cutoff operator $H_{\Omega,n}$ is conservative on $L_\infty(\Omega_n)$  by Step~1.

\smallskip

\noindent{\bf Step~3}$\;$
The third and fourth steps consist of removing the cutoff by a suitable limit $n\to\infty$,
first by  $L_2$-arguments  and then by  $L_1$-arguments.

It is convenient to view $H_\Omega$ and $H_{\Omega, n}$ as symmetric operators on $L_2(\Ri^d)$.
Since  the coefficients
$c_{ij}$ of $H_\Omega$ are in $W^{1,\infty}(\Omega)$ the operator 
  can be extended  to a symmetric operator on the domain $C_c^\infty(\Ri^d)$.
The extension corresponds to the operator $H_\Omega\oplus \,0$ with domain $C_c^\infty(\Omega)\oplus L_2(\Omega^{\rm c})$.
The Markov uniqueness of $H_\Omega$ on $L_2(\Omega)$ implies that the extended operator has a unique submarkovian 
extension $H^F_\Omega\oplus \,0$ on $L_2(\Ri^d)$ and for simplicity of notation  we set $H=H^F_\Omega\oplus\, 0$.
Similarly, since $H_{\Omega, n}$ is Markov unique by Step~2 there is a unique submarkovian operator $H_n=H^F_{\Omega,n}\oplus \,0$
which extends $H_{\Omega,n}$.
We let $h$ and $h_n$ denote the corresponding Dirichlet forms on $L_2(\Ri^d)$.

The $\rho_n$ form an increasing sequence of functions on $\Ri^d$, by definition.
Therefore  the~$h_n$ are a monotonically  increasing family of forms on $L_2(\Ri^d)$.
This implicitly uses the Markov uniqueness through the identification $l_\Omega=\overline h_\Omega$ and hence $l_{\Omega, n}=\overline h_{\Omega,n}$.
Therefore one can define $h_\infty$ by $D(h_\infty)=\bigcap_{n\geq1}D(h_n)$ and $h_\infty(\varphi)=\sup_{n\geq1}h_n(\varphi)$
for all $\varphi\in D(h_\infty)$.
The form $h_\infty$ is closed (see, for example, \cite{Kat1}, Section~VIII.3.4)
and since the $h_n$ are Dirichlet forms the supremum $h_\infty$ is also   a Dirichlet form.
Moreover,  by direct calculation $h_\infty(\varphi)=h(\varphi)$ for all $\varphi\in C_c^\infty(\Ri^d)$.
Hence $h_\infty\supseteq h$.
Then it follows from the monotone convergence of the forms $h_n$ that one has strong $L_2$-convergence of the resolvents
$(\lambda I+H_n)^{-1}$ to the resolvent $(\lambda I+H_\infty)^{-1}$ for all $\lambda>0$
where $H_\infty$ is the submarkovian operator corresponding to the form $h_\infty$.
Hence 
\[
H_\infty(I+\varepsilon H_\infty)^{-1}\varphi=\lim_{n\to\infty}H_n(I+\varepsilon H_n)^{-1}\varphi=\lim_{n\to\infty}(I+\varepsilon H_n)^{-1}H\varphi
\]
for all $\varphi\in C_c^\infty(\Ri^d)$.
Since $(I+\varepsilon H_n)^{-1}$ converges strongly to the identity operator as $\varepsilon\to0$ it follows that $C_c^\infty(\Ri^d)\subseteq D(H_\infty)$ 
and $H_\infty\varphi=H\varphi$
for all $\varphi\in C_c^\infty(\Ri^d)$.
Thus $H_\infty$ is a submarkovian extension of $H$ and by Markov uniqueness one has $H_\infty=H$.

The foregoing arguments establish that the $H_n$ converge to $H$ in the strong resolvent sense on $L_2(\Ri^d)$.
Therefore the submarkovian semigroups $S^{(n)}_t$ generated by the $H_n$ converge strongly on $L_2(\Ri^d)$ to the submarkovian semigroup $S_t$
generated by $H$.

Note that by construction the semigroup $S^{(n)}_t$ leaves both $L_2(\Omega_n)$  and the orthogonal complement $L_2(\Omega^{\rm c})$ invariant.  
The semigroup is conservative on $L_\infty(\Omega)$ by Step~2 and is equal to the identity semigroup on the orthogonal complement.
Therefore the $S^{(n)}_t$ are conservative semigroups on $L_\infty(\Ri^d)$ which are strongly $L_2$-convergent to $S_t$.
But this is not sufficient to ensure that $S_t$ is conservative.
For this one needs $L_1$-convergence.

\smallskip

\noindent{\bf Step~4}$\;$
The fourth step in the proof consists in proving that the semigroups $S^{(n)}_t$ are strongly convergent on $L_1(\Ri^d)$ to 
$S_t$ (see \cite{RSi2}, Proposition~6.2, for a similar result).

 Since the semigroups $S^{(n)}_t$ and $S_t$  are all submarkovian  it suffices to prove convergence
 on a  subset of $L_1$ whose span is dense.
 In particular it suffices to prove convergence on $L_1(A)\cap L_2(A)$ for each bounded open
 subset $A$ of $\Omega$.
 Moreover one can restrict to positive functions.

 Fix $A\subset\Omega$ and $\varphi_A\in L_1(A)\cap L_2(A)$.
 Assume $\varphi_A$ is positive.
 Next let $B\supset A$ be  a bounded closed set of $\Ri^d$.
Then
\begin{eqnarray}
\|(S^{(n)}_t-S_t)\varphi_A\|_1&\leq &
\|\one_B(S^{(n)}_t-S_t)\varphi_A\|_1
+\|\one_{B^{\rm c}}S^{(n)}_t\varphi_A\|_1+\|\one_{B^{\rm c}}S_t\varphi_A\|_1\nonumber\\[5pt]
&\leq&|B|^{1/2}\|(S^{(n)}_t-S_t)\varphi_A\|_2
+|(\one_{ B^{\rm c}},S^{(n)}_t\varphi_A)|+|(\one_{ B^{\rm c}},S_t\varphi_A)|
\label{esm5.11}
\end{eqnarray}
where we have used the positivity of the semigroups and the functions to express the $L_1$-norms as pairings between $L_1$ and $L_\infty$.
Therefore it suffices to prove that the last two terms can be made arbitrarily small, uniformly in $n$, by suitable choice of $B$.
Then the $L_1$-convergence follows from the $L_2$-convergence of Step~3.
But the uniform estimate follows by Davies--Gaffney bounds using the arguments of Proposition~3.6 of \cite{ERSZ1}.
We briefly sketch the proof.

First  one can  associate a set theoretic (quasi-)distance with a quite general Dirichlet form (see, for example, \cite{BM} \cite{Stu2} \cite{AH} or \cite{ERSZ2}).
Specifically  we introduce the family of 
 $d_n(X\,;Y)$  of distances corresponding to the Dirichlet forms $h_n$
and a similar distance $d(X\,;Y)$ corresponding to  $h$ following the definitions of Section~1 of \cite{ERSZ2}.
Here $X$ and $Y$ are measurable subsets of $\Ri^d$ and  $d_n(X\,;Y)\in [0,\infty]$.
The definition of the distances is quite technical and $d_n(X\,;Y)$  takes the value $+\infty$ if $X$ or $Y$ is not a subset of $\Omega_n$.
But since $h_n\leq h$ one has $d_n(X\,;Y)\geq d(X\,;Y)$ and since $h(\varphi)\leq \|C\|\|\nabla\varphi\|_2^2$ one also has
\[
 d(X\,;Y)\geq  \|C\|^{-1/2}|X-Y|= \|C\|^{-1/2}\inf_{x\in X, \,y\in Y}|x-y|
 \]
 (see \cite{ERSZ2}, Section~5, for a discussion of the monotonicity properties of the distances).
 Then the Davies--Gaffney bounds \cite{Gaf} \cite{Dav12} \cite{Stu2} \cite{AH}  as presented  in Theorem~2 of \cite{ERSZ2} give
 \begin{eqnarray}
|(\varphi_{X},S^{(n)}_t\varphi_{Y})|&\leq& e^{-d_n(X ;Y)^2 (4t)^{-1}}\|\varphi_{X}\|_2\|\varphi_{Y}\|_2\nonumber\\
&\leq& e^{-d(X ;Y)^2 (4t)^{-1}}\|\varphi_{X}\|_2\|\varphi_{Y}\|_2\leq e^{-|X -Y|^2 (4\|C\|t)^{-1}}\|\varphi_{X}\|_2\|\varphi_{Y}\|_2
\label{esm5.12}
\end{eqnarray}
for all  $\varphi_{X}\in L_2(X)$ and $\varphi_{Y}\in L_2(Y)$.
These bounds are uniform in $n$ and are conveniently expressed in terms of the Euclidean distance.

Now choose $R$ sufficiently large that $A\subseteq B_{R}=\{x:|x|<R\}$ and let 
$B=\overline{B_{2R}}$. 
Then one can separate $B^{\rm c}$ into annuli $B_{(n+1)R}\backslash B_{nR}$ and make a quadrature estimate, as in the proof of Proposition~3.6 of \cite{ERSZ1}, to find
\begin{eqnarray*}
e^{-|A -B^{\rm c}|^2 (4\|C\|t)^{-1}} \leq & \sum_{n\geq2}e^{-|B_{R} - B_{nR}^{\rm c}|^2 (4\|C\|t)^{-1}}|B_{(n+1)R}|^{1/2} 
\leq a\,R^{d/2}e^{-bR^2t^{-1}}
  \end{eqnarray*}
with $a, b>0$.
Therefore combining these bounds with (\ref{esm5.11}) and (\ref{esm5.12}) one obtains the equicontinuous bounds
\[
\|(S^{(n)}_t-S_t)\varphi_A\|_1\leq a'\,R^{d/2}\,\|(S^{(n)}_t-S_t)\varphi_A\|_2+2\,a\,R^{d/2}e^{-bR^2t^{-1}}
\]
where $a',a,b>0$ are all independent of $n$.
It follows immediately that $\|(S^{(n)}_t-S_t)\varphi_A\|_1\to0$ as $n\to\infty$.
Thus the $S^{(n)}_t$ converge strongly to $S_t$ on $L_1(\Ri^d)$ and in particular on the invariant  subspace $L_1(\Omega)$.

\smallskip

\noindent{\bf Step~5}$\;$ Finally we combine the conclusions of Steps~1 and 4 to deduce that $S^F_t$ is conservative on $L_\infty(\Omega)$.

It follows from Step~1 that the  semigroup generated by the Friedrichs extension $H^F_{\Omega,n}$ of the cutoff operator $H_{\Omega,n}$  is a conservative semigroup on $L_\infty(\Omega_n)$.
Therefore the extension $S^{(n)}_t$ of the semigroup  to $L_\infty(\Ri^d)$ is also conservative since
\[
S^{(n)}_t\one=(S^{(n)}_t\oplus I)(\one_{\Omega_n}\oplus\one_{\Omega_n^{\rm c}})=S^{(n)}_t\one_{\Omega_n}\oplus\one_{\Omega_n^{\rm c}}=
\one_{\Omega_n}\oplus\one_{\Omega_n^{\rm c}}=\one
\;.
\]
Then, however,
\[
(\one,\varphi)=\lim_{n\to\infty}(\one, S^{(n)}_t\varphi)=(\one,S_t\varphi)
\]
for all $\varphi\in L_1(\Ri^d)$ by Step~4.
Hence $S_t$ is conservative on $L_\infty(\Ri^d)$ and its restriction $S^F_t$ to the invariant subspace $L_\infty(\Omega)$ is conservative.

\medskip

\noindent
\ref{tsm1.2-0}$\Leftrightarrow$\ref{tsm1.2-2}$\;$ This follows by an argument of Davies, \cite{Dav14} Theorem~2.2,
which was given for operators with smooth coefficients but which is also valid for operators
with Lipschitz coefficients.
In fact Davies argues that $S^F_t$ is conservative if and only if $C_c^\infty(\Omega)$ is a core for the generator of the semigroup acting on $L_1(\Omega)$.
But this is equivalent to $L_1$-uniqueness (see \cite{Ebe}, Section~1b).
Davies arguments need a slight modification to cover the operator $H_\Omega$ but this is not difficult by the discussion of 
elliptic regularity properties in  Section~\ref{S3}.
We omit further details.

\smallskip

\noindent
\ref{tsm1.2-2}$\Rightarrow$\ref{tsm1.2-3}$\;$ This is a general feature which is proved in \cite{Ebe}, Lemma~1.6.
\hfill$\Box$

\bigskip

Note that the implication \ref{tsm1.2-0}$\Rightarrow$\ref{tsm1.2-3}, which is an indirect consequence of the foregoing proof, 
can be easily deduced from Theorem~\ref{tsm1.0}.
Let $S_t$ denote the semigroup generated by $L_\Omega$.
Then $S^F_t\varphi\leq S_t\varphi$ for all positive $\varphi\in L_2(\Omega)$ and all $t\geq0$ by Theorem~\ref{tsm1.0}.\ref{tsm1.0-5}.
But if $\varphi\in L_2(\Omega)$ and $0\leq \varphi\leq \one_\Omega$ then 
\[
\one_\Omega=S^F_t\one_\Omega=S^F_t\varphi+S^F_t(\one_\Omega-\varphi)\leq
S_t\varphi+S_t(\one_\Omega-\varphi)=S_t\one_\Omega\leq\one_\Omega\;.
\]
Therefore the inequalities are equalities 
 and $S^F_t\varphi=S_t\varphi$ for all positive $\varphi\in L_2(\Omega)$ such that  $0\leq \varphi\leq \one_\Omega$ and for all $t\geq0$. 
 It follows immediately that $S^F_t=S_t$ for all $t\geq0$.
Therefore $H_\Omega^F=L_\Omega$ and  $H_\Omega$ is Markov unique by Theorem~\ref{tsm1.0}.\ref{tsm1.0-2}.
(This argument follows  the latter part of the  proof of Corollary~3.4 in \cite{Ebe}.)

\section{Concluding remarks}\label{S6}

In this concluding section we discuss various results and examples concerning $L_p$-uniqueness, sets of capacity zero and irreducibility properties.

\subsection{$L_p$-uniqueness}

First note that Lemma~\ref{lsm5.1} gives a condition,
in terms of an approximation to the identity, which ensures that  $S^F_t$ is conservative, and consequently  $H_\Omega$ is  $L_1$-unique.
But if $p\in[1,2]$ there is a similar sufficient condition  for $L_p$-uniqueness.

\begin{prop}\label{psm6.1}
Assume $p\in[1,2\,]$.
If there exists a sequence $\eta_n\in C_c^\infty(\Omega)$
such that  $0\leq \eta_n\leq \one_\Omega$, $\|(\eta_n-\one_{\Omega})\psi\|_2 \to0$ for all $\psi\in L_2(\Omega)$   and $\|\Gamma(\eta_n)\|_{p/(2-p)}\to0$ as $n\to\infty$
then $H_\Omega$ is $L_p$-unique.
\end{prop}

In the case $p=2$   Davies has established similar criteria (see \cite{Dav14}, Theorems~3.1 and~3.2).
(If $p=2$ then $p/(2-p)$ is understood to be $\infty$.)
Moreover,  if $p=1$ then $\|\Gamma(\eta_n)\|_1=h_\Omega(\eta_n)$ and the condition for $L_1$-uniqueness agrees with the condition in Lemma~\ref{lsm5.1}.

Proposition~\ref{psm6.1}  is essentially a corollary of the following.

\begin{lemma}\label{tse3.2}
If  $\varphi\in D(H_\Omega^*)$ and $\eta\in  C_c^\infty(\Omega)$
then
 \begin{equation}
(\eta^2\varphi,H_\Omega^*\varphi)\geq -(\varphi,\Gamma(\eta)\varphi)
 \label{ebd4.1}
 \end{equation}
 where $\Gamma$ is the {\it carr\'e du champ} associated with $H_\Omega$.
 Therefore if $(I+H_\Omega^*)\varphi=0$  then
 \begin{equation}
 \|\eta\, \varphi\|_2^2\leq (\varphi,\Gamma(\eta)\varphi)
\label{ebd4.2}
 \end{equation}
for all $\eta\in  C_c^\infty(\Omega)$.
 \end{lemma}
 \proof\
First, if $\eta\in C_c^\infty(\Omega)$ and $\varphi\in D(H_\Omega^*)$ then
 $\eta\varphi , \eta^2\varphi \in  D(\overline H_\Omega)$
by  Theorem~\ref{tse3.1}.
Therefore
\begin{eqnarray*}
 2\,(\eta^2\varphi,H^*_\Omega \varphi)&=&2\,\RRe(\eta^2\varphi,H^*_\Omega \varphi)\\
 &=&(\overline H_\Omega \eta^2\varphi,\varphi)+(\varphi,\overline H_\Omega \eta^2\varphi)\\
&\geq& (\overline H_\Omega \eta^2\varphi,\varphi)+(\varphi,\overline H_\Omega\eta^2\varphi)-2\,(\overline H_\Omega\eta\varphi,\eta\varphi) 
\end{eqnarray*}
since $\overline H_\Omega\geq 0$.
But if $\Omega'\subset\subset\Omega$ is bounded and $\supp\eta\subset\Omega'$ then one may construct the strongly elliptic extension $L$ of $H_{\Omega'}$ to $L_2(\Ri^d)$
as in  the  proof of Theorem~\ref{tse3.1}.
Then since $\overline H_\Omega\eta\varphi=\overline H_{\Omega'}\eta\varphi$ etc. one has 
\begin{eqnarray*}
(\overline H_\Omega \eta^2\varphi,\varphi)+(\varphi,\overline H_\Omega\eta^2\varphi)-2\,(\overline H_\Omega\eta\varphi,\eta\varphi) 
&=&(\varphi, \eta^2 L(\varphi))+(\varphi,L(\eta^2\varphi))-2\,(\varphi,\eta L(\eta\varphi))\\
&=&(\varphi, L(\eta^2)\varphi)-2\,(\varphi, \eta L(\eta)\varphi)=-2\,(\varphi,\Gamma(\eta)\varphi)
\end{eqnarray*}
where  we have used the distributional relation (\ref{ereg2.1}) several times.
Combination 
 of the last two estimates  immediately yields
 (\ref{ebd4.1}).
\hfill$\Box$

\begin{remarkn}\label{rbd3.1}
The essence of the foregoing calculation is the formal double commutator identity
\[
(\ad\,\eta)^2(H_\Omega)=[\eta,[\eta,H_\Omega]]=-2\,\Gamma(\eta)\;.
\]
Double commutator estimates of a different nature were used to prove general
self-adjointness results in \cite{Rob7}, e.g.\ Theorem~2.10, (see also \cite{ER27}, Proposition~2.3).
\end{remarkn}

\noindent{\bf Proof of Proposition~\ref{psm6.1}}$\;$
It suffices to prove that the range of $I+H_\Omega$ is dense in $L_p(\Omega)$.
Therefore assume that $\varphi\in L_q(\Omega)$, the dual space of $L_p(\Omega)$, and $(I+H_\Omega^*)\varphi=0$.
Since $q\in[2,\infty]$ it follows that $\eta_n\varphi=-\eta_n H^*_\Omega\varphi\in L_2(\Omega)$ and then
 (\ref{ebd4.1}) gives
\[
 \|\eta_n \varphi\|_2^2\leq (\varphi,\Gamma(\eta_n)\varphi)=\int\Gamma(\eta_n)\,\varphi^2\leq \|\varphi\|_q^2\,\|\Gamma(\eta_n)\|_{p/(2-p)}
 \;.
 \]
Taking the limit $n\to \infty$ one deduces that $\|\varphi\|_2=0$ so $\varphi=0$ and the range is dense.
\hfill$\Box$

\bigskip

If $p=2$  then the statement of Proposition~\ref{psm6.1} be strengthened.

\begin{cor}\label{csm6.1}
Assume there exists a sequence $\eta_n\in C_c^\infty(\Omega)$
such that $0\leq \eta_n\leq \one_\Omega$, $\|(\eta_n-\one_{\Omega})\psi\|_2 \to0$ for all $\psi\in L_2(\Omega)$   and $\sup_{n\geq1}\|\Gamma(\eta_n)\|_\infty<\infty$. 
Then $H_\Omega$ is $L_2$-unique, i.e.\ $H_\Omega$ is essentially self-adjoint.
\end{cor}
\proof\
It suffices to prove that the range of $I+\varepsilon H_\Omega$ is dense in $L_2(\Omega)$ for all small $\varepsilon>0$.
But if $\varphi\in D(H^*_\Omega)$ and $(I+\varepsilon  H_\Omega^*)\varphi=0$ then the foregoing argument gives 
\[
 \|\eta_n \varphi\|_2^2\leq \varepsilon\,(\varphi,\Gamma(\eta_n)\varphi)\leq\varepsilon \sup_{n\geq1}\|\Gamma(\eta_n)\|_\infty\,\|\varphi\|_2^2
 \;.
 \]
Therefore $\|\varphi\|_2=0$ for all small $\varepsilon>0$.
Thus $\ker(I+\varepsilon\,H_\Omega)=\{0\}$ and the range of $I+\varepsilon  H_\Omega$ is dense.
\hfill$\Box$

\begin{exam}\label{exsm5.1}
Let $\Omega=\Ri^d$.
Then the operator $H=-\sum^d_{i,j=1}\partial_ic_{ij}\partial_j$ acting on $C_c^\infty(\Ri^d)$ with $c_{ij}\in W^{1,\infty}(\Ri^d)$ and $(c_{ij})>0$ is $L_2$-unique
as a consequence of  Proposition~\ref{psm6.1}.
It suffices to choose $\eta_n\in C_c^\infty(\Ri^d)$ with $0\leq\eta_n\leq 1$, $\eta_n(x)=1$ if $|x|\leq n$ and $\|\nabla\eta_n\|_\infty\leq a\,n^{-1}$.
Then the $\eta_n$ converge pointwise to the identity as $n\to\infty$ and $\|\Gamma(\eta_n)\|_\infty\leq a\,n^{-2}\|C\|\to0$.
The $L_2$-uniqueness  implies that $H$ is Markov unique. 
Therefore  $H$ is also $L_1$-unique and $S^F_t$ is conservative by Theorem~\ref{tsm1.2}.
\end{exam}

\begin{exam}\label{exsm5.2}
Assume that $c_{ij}\in W^{1,\infty}(\Omega)$ and $0<C(x)\leq a\,d_{\partial\Omega}(x)^2$ for some $a>0$ and all $x\in \Omega$ where $d_{\partial\Omega}$ is the Euclidean
distance to the boundary $\partial\Omega$ of $\Omega$.
Then $H_\Omega$ is $L_2$-unique.
Again this follows from Proposition~\ref{psm6.1}.
It suffices to define $\eta_n$ as in the proof of Proposition~\ref{psm4.1}.
Then the $\eta_n$ converge pointwise to $\one_\Omega$ and 
$\|\Gamma(\eta_n)\|_\infty\leq a\,(\log n)^{-2}$.
Therefore $L_2$-uniqueness follows from Proposition~\ref{psm6.1} with $p=2$.
More generally if  $d(\partial\Omega)=d-1$  one can 
use the calculational procedure of the proof of
Proposition~\ref{psm4.1} to deduce that if $p\in[1,2\,]$ and $0<C(x)\leq a\,d_{\partial\Omega}(x)^{(3p-2)/p}$ then $H_\Omega$ is $L_p$-unique.
In particular $L_1$-uniqueness follows if $0<C(x)\leq a\,d_{\partial\Omega}(x)$.
\end{exam}

Although  the approximation criteria  for $L_1$-uniqueness  and $L_2$-uniqueness in Proposition~\ref{psm6.1}
are superficially similar they are of a totally different geometric character.
The first involves the  norm $\|\Gamma(\eta)\|_1$ which is related
to the capacity and the second  involves  the norm $\|\Gamma(\eta)\|_\infty$ which  is related to the Riemannian distance.
In one-dimension the first estimate is optimal but the second is suboptimal.
This is illustrated by the following   example adapted from \cite{CMP}  (see also \cite{Ebe} \cite{RSi3}).

\begin{exam}\label{exsm5.3}
Assume $d=1$ and $\Omega=\langle-1,1\rangle$. 
Further let  $H$ be the operator with domain $C_c^\infty(-1,1)$ and action $H\varphi=-(c\,\varphi')'$
where $c(x)=(1-x^2)^\delta$.
Then $c\in W^{1,\infty}(-1,1)$ if and only if $\delta\geq1$.
Set $W(x)=\int^x_0c^{-1}$.
Thus $H^*W=0$.
It follows that  $H$ is $L_p$-unique for $p\in[1,\infty\rangle$ if and only if $W\not\in L_q(-1,1)$ where $q$ is conjugate to $p$ (see \cite{CMP} Proposition~3.5).
Hence $H$ is $L_1$-unique for all $\delta\geq 1$ and $L_p$-unique for $p>1$  if and only if $\delta>(2p-1)/p$.
In particular it is $L_2$-unique  if and only if $\delta>3/2$ and  $L_p$-unique for all $p\in[1,\infty\rangle$ if and only if $\delta\geq2$.
Alternatively, $H$ is Markov unique for all $\delta\geq1$ by \cite{Ebe}, Theorem~3.5.
Thus Markov uniqueness and $L_1$-uniqueness are simultaneously valid in agreement with Theorem~\ref{tsm1.2}.

The $L_1$-uniqueness can be verified by the criterion of Proposition~\ref{psm6.1}.
Define $\eta_n$ by $\eta_n(x)=1-W(x)/W(1-n^{-1})$ if $x\in[0,n^{-1}\rangle$, $\eta_n(x)=0$ if $x\geq 1-n^{-1}$ and $\eta_n(-x)=\eta_n(x)$ for all $x\geq0$.
Since $\delta\geq1$ it follows that $\eta_n$ converges monotonically upward to $\one_{\langle-1,1\rangle}$ as $n\to\infty$.
But $\Gamma(\eta_n)=c\,|\eta_n'|^2$.
Thus $h(\eta_n)=\|\Gamma(\eta_n)\|_1=2\,W(1-n^{-1})^{-1}\to0$ as $n\to\infty$.
Therefore   $L_1$-uniqueness of $H$ follows for all $\delta\geq1$. 
But  $\|\Gamma(\eta_n)\|_\infty\sim n^{(2-\delta)}$ and this is bounded if and only if $\delta\geq2$.
Therefore the $L_2$-uniqueness only follows for $\delta\geq 2$ and not for the full range $\delta>3/2$.

Note that the Riemannian distance corresponding to the metric $c^{-1}$ is given by $d(x\,;y)=|\int^x_yc^{-1/2}|$.
Thus the distance from the origin to the boundary, $d(0\,;1)=d(0\,;-1)$,  is finite for all $\delta\in[1,2\rangle$.
Therefore if $\delta\in\langle 3/2, 2\rangle$ then the distance to the boundary is finite but $H$ is nonetheless essentially self-adjoint.
\end{exam}

\subsection{Sets of capacity zero}

Let $A$ be a closed subset of $\overline\Omega$ with $|A|=0$.
In this subsection we  assume that the coefficients $c_{ij}$
are real, symmetric, $c_{ij}\in W^{1,\infty}(\Omega)$ and $C(x)>0$ for all $x\in\Omega\backslash A$.
Then we define the operators $H_\Omega$ and $H_{\Omega\backslash A}$ with the coefficients $c_{ij}$
on $C_c^\infty(\Omega)$ and $C_c^\infty(\Omega\backslash A)$, respectively.
All the foregoing considerations apply to $H_{\Omega\backslash A}$  because the matrix of coefficients $C$
is non-degenerate on $\Omega\backslash A$ but they do not necessarily apply to $H_\Omega$ since $C$
can be degenerate on $A$.
Nevertheless $H_\Omega\supseteq H_{\Omega\backslash A}$.
Hence uniqueness criteria for $H_{\Omega\backslash A}$ give sufficient conditions
for uniqueness of $H_\Omega$.
For example if $H_{\Omega\backslash A}$ is Markov unique then $H_\Omega$ is Markov unique.
But Markov uniqueness of $H_{\Omega\backslash A}$ is equivalent to the boundary $\partial(\Omega\backslash A)$
having zero capacity and this is equivalent to $\partial\Omega$ and $A$ both having zero capacity.
Thus the boundary condition $\capp_\Omega(\partial\Omega)=$ is sufficient for $H_\Omega$ to be Markov unique
if in addition the degeneracy set $A$ has zero capacity.
This typically occurs for one of two reasons.
Either $d(A)\leq d-2$ and $\capp_\Omega(A)=0$ independently of the behaviour of the coefficients  in the neighbourhood
of $A$ or $d(A)$ is arbitrary and the coefficients have a correspondingly strong degeneracy on $A$ (see Proposition~\ref{psm4.1}).
We illustrate these possibilities with two simple examples.

\begin{exam}\label{exsm5.4}
Let $\Omega=\Ri^d$ and consider the operator 
$H=-\sum^d_{i,j=1}\partial_ic_{ij}\partial_j$ acting on $C_c^\infty(\Ri^d)$ with $c_{ij}\in W^{1,\infty}(\Ri^d)$ and $(c_{ij})>0$
on the complement $A^{\rm c}$ of a closed set $A$ with $d(A)\leq d-2$.
Then $H$ is $L_1$-unique and Markov unique.
This follows because $\partial A^{\rm c}=A$ and  $\capp_\Omega(A)=0$ by the estimates of Proposition~\ref{psm4.1}.
Therefore $H$ is $L_1$-unique and Markov unique on $C_c^\infty(\Ri^d\backslash A)$ by Theorems~\ref{tsm1.1} and \ref{tsm1.2}.
\end{exam}

\begin{exam}\label{exsm5.5}
Again let $\Omega=\Ri^d$ and 
$H=-\sum^d_{i,j=1}\partial_ic_{ij}\partial_j$ the operator acting on $C_c^\infty(\Ri^d)$ with coefficients  $c_{ij}\in W^{1,\infty}(\Ri^d)$ and $(c_{ij})>0$
on the complement $A^{\rm c}=\Ri^d\backslash A$ of a closed set $A$ with the property that $A^{\rm c}=\Omega_1\cup \Omega_2$ with $\Omega_1\cap\,\Omega_2=\emptyset$.
Now assume $\capp_\Omega(A)=0$.
Since $A=\partial(\Omega_1\cup\Omega_2)=\overline\Omega_1\cap\,\overline\Omega_2$
it follows that 
 $H$ is $L_1$-unique and Markov unique on $C_c^\infty(\Omega_1\cup\,\Omega_2)$.
 Therefore $H$ is $L_1$-unique and Markov unique on $C_c^\infty(\Ri^d)$.
 Moreover, the unique Markov extension, the Friedrichs extension $H^F$, must coincide with the Friedrichs extension of $H$ on
  $C_c^\infty(\Omega_1\cup\,\Omega_2)=C_c^\infty(\Omega_1)+C_c^\infty(\Omega_2)$.
  But it follows readily from the definition of the Friedrichs extension  
that this latter operator is of the form $H^F_{\Omega_1}\oplus H^F_{\Omega_2}$ on $L_2(\Omega_1)\oplus L_2(\varphi_2)$.
Therefore  the semigroup $S^F_t$ generated by $H_F$ leaves the subspaces $L_2(\Omega_1)$ and $L_2(\Omega_2)$ invariant.
\end{exam}

\subsection{Irreducibility and ergodicity}

In Example~\ref{exsm5.4} the set $\Ri^d\backslash A$ on which the coefficients of the operator $H$
are non-degenerate has two disjoint components $\Omega_1$ and $\Omega_2$.
Consequently the corresponding Markov semigroup has two invariant subspaces $L_2(\Omega_1)$ and $L_2(\Omega_2)$.
We conclude by giving a general result that relates connectedness of the set of non-degeneracy and ergodicity of the 
corresponding Friedrichs semigroup.

The absence of non-trivial invariant subspaces is variously defined as ergodicity or irreducibility of a semigroup.
The property  can be characterized by strict positivity.
In particular the positive semigroup $S_t$ on $L_2(\Omega)$  is defined to be irreducible
if for every $t>0$ and every positive, nonzero, $ \varphi\in L_2(\Omega)$ one has $S_t\varphi>0$ almost everywhere
(see, for example, \cite{Ouh5}, Definition~2.8).
This is clearly  equivalent to the requirement that 
 $(\varphi, S_t\psi)>0$ for all positive, nonzero, $ \varphi, \psi\in L_2(\Omega)$ and
for all $t>0$.

Now consider the submarkovian semigroups generated by extensions of $H_\Omega$ always under the assumptions
of Theorem~\ref{tsm1.0}.
The following proposition extends Theorem~4.5 of \cite{Ouh5} to this situation.

\begin{prop}\label{psm5.1}
Let $S^F_t$ denote the semigroup generated by the Friedrichs extension $H^F_\Omega$ of $H_\Omega$.

The following conditions are equivalent:
\begin{tabel}
\item\label{psm5.1-1}
$S^F_t$ is irreducible,
\item\label{psm5.1-2}
$\Omega$ is connected.
\end{tabel}
Moreover if these conditions are satisfied then each semigroup generated by  a submarkovian extension of $H_\Omega$ is irreducible.
\end{prop}
\proof\
First note that if $K_\Omega$ is a submarkovian extension of $H_\Omega$ then the semigroup $e^{-tK_\Omega}$ dominates $S^F_t$ by Theorem~\ref{tsm1.0}.\ref{tsm1.0-5}.
It follows immediately that irreducibility of $S^F_t$ implies irreducibility of $e^{-tK_\Omega}$.

\smallskip
\noindent\ref{psm5.1-1}$\Rightarrow$\ref{psm5.1-2}$\;$
This follows by the foregoing discussion. 
If $\Omega$ is not connected then $S^F_t$ is not irreducible.

\smallskip

\noindent\ref{psm5.1-2}$\Rightarrow$\ref{psm5.1-1}$\;$
Let $\Omega_n$ be an increasing family of connected open subsets of $\Omega$ with $\Omega_n\subset\subset\Omega_{n+1}$
and $\Omega=\bigcup_{n\geq1}\Omega_n$.
Then $H_n=H_\Omega|_{C_c^\infty(\Omega_n)}$ is a strongly elliptic operator on $L_2(\Omega_n)$.
Let $H_n^F$ denote the Friedrichs extension of $H_n$ and $S^{(n)}_t$ the corresponding semigroup.
The extension $H_n^F$  corresponds to  Dirichlet boundary conditions on $H_n$.
It follows that $S^{(n)}_t$ is irreducible  on $L_2(\Omega_n)$ by \cite{Ouh5}, Theorem~4.5.
But $S^F_t\varphi\geq S^{(n+1)}_t\one_{\Omega_{n+1}}\varphi \geq S^{(n)}_t\one_{\Omega_{n}}\varphi$
for all positive $\varphi\in L_2(\Omega)$ by Corollary~2.3 of \cite{ER30}.
Therefore $S^F_t$ is irreducible on $L_2(\Omega)$.~\hfill$\Box$

\bigskip

\subsection*{Acknowledgement}
Our interest in the characterization of operator extensions  by their boundary behaviour  originated in a series of interesting discussions
 with Ricardo Weder and Gian Michele Graf whilst the first author was visiting 
the Institut f\"ur Theoretische Physik at the Eidgen\"ossische Technische Hochschule, Z\"urich in 2007.
The author is indebted to J\"urg Fr\"ohlich and Gian Michele Graf  for facilitating this visit and to the ETH for providing financial support.
Both authors are grateful to Michael R\"ockner and Masayoshi Takeda for comments on an earlier draft of the paper.

\appendix
\section{Strong ellipticity}\label{S2}

In this appendix we  recall some basic properties of the operator $L=-\sum^d_{i,j=1}\partial_i\,c_{ij}\,\partial_j$,
with real symmetric Lipschitz continuous  coefficients  $c_{ij}$  and domain $D(L)=C_c^\infty(\Ri^d)$, acting on $L_2(\Ri^d)$
under the  hypothesis of strong ellipticity.

\begin{prop}\label{pse1}$\;$Assume $c_{ij}=c_{ji}\in W^{1,\infty}(\Ri^d)$ and  that  $C=(c_{ij})\geq \mu I>0$ uniformly over $\Ri^d$.
Then one has the following.
\begin{tabel}
\item\label{pse1-1}
$L$  is essentially self-adjoint,
\item\label{pse1-2}
$D(\overline L)=W^{2,2}(\Ri^d)$,
\item\label{pse1-3}
$(I+\overline L)^{-1}$ is a bounded operator from $W^{-\delta,2}(\Ri^d)$ to $W^{2-\delta,2}(\Ri^d)$ for all $\delta\in[0,1]$.
\end{tabel}
\end{prop}
 The conclusions of the proposition are well known.
If the operator is strongly elliptic and $c_{ij}\in W^{2,\infty}(\Ri^d)$ then the result  follows in its entirety  from  \cite{Rob7} or \cite{ER27}
but we have not
found a suitable reference for the complete statement with $c_{ij}\in W^{1,\infty}(\Ri^d)$.
We  briefly sketch the proof.

\medskip

\noindent{\bf Sketch of proof of Proposition~\ref{pse1}}$\;$ First,  since the operator $L$  is symmetric on $L_2(\Ri^d)$ it is closable
and  its closure $\overline L$  is self-adjoint if and only if 
 the range condition $R(\kappa I+\overline L)=L_2(\Ri^d)$ is satisfied  for  large positive $\kappa$.
Since the coefficients $c_{ij}\in W^{1,\infty}(\Ri^d)$ the latter condition can be established by 
the following variant of Levi's  parametrix argument (see, for example, \cite{Fri}).

Fix $y\in\Ri^d$ and introduce the constant coefficient operator $L_y=-\sum^d_{i,j=1}c_{ij}(y)\,\partial_i\partial_j$
with domain $W^{2,2}(\Ri^d)$.
Then $L_y$ is a positive self-adjoint operator which generates a translationally invariant submarkovian semigroup
with an integral kernel 
\begin{equation}
K^{(y)}_t(x)=(\det C(y))^{-1}(4\pi t)^{-d/2}e^{-(x,C(y)^{-1}x)/4t}
\;.
\label{ese1.40}
\end{equation}
The kernel $R^{(y)}_\kappa$ of the resolvent $(\kappa I+L_y)^{-1}$ is then given by
\[
R^{(y)}_\kappa(x)=\int^\infty_0dt\,e^{-\kappa t}K^{(y)}_t(x)
\;.
\]
Next define $R_\kappa$ as a bounded   operator on the spaces $L_p(\Ri^d)$ by 
\[
(R_\kappa\varphi)(x)=\int_{\Ri^d}dy\, R^{(y)}_\kappa(x-y)\,\varphi(y)
\;.
\]
It follows  that 
if $\kappa\geq 1$ then 
$R_\kappa L_2(\Ri^d)\subseteq D(\overline L)$ and 
\[
(\kappa I+\overline L)R_\kappa=I+Q_\kappa
\]
 where the  $Q_\kappa$ are the  bounded operators 
 \begin{eqnarray*}
 (Q_\kappa\varphi)(x)&=&-\sum^d_{i,j=1}\int_{\Ri^d}dy\,\Big(\partial_i(c_{ij}(x)-c_{ij}(y))\partial_j R^{(y)}_\kappa\Big)(x-y)\,\varphi(y)\\
 &&\hspace{-1.2cm}=-\sum^d_{i,j=1}\int_{\Ri^d}dy\,\Big((\partial_ic_{ij})(x)(\partial_jR^{(y)}_\kappa)(x-y)+(c_{ij}(x)-c_{ij}(y))(\partial_i\partial_j R^{(y)}_\kappa)(x-y)\Big)\,\varphi(y)
 \;.
 \end{eqnarray*}
 Since the coefficients $c_{ij}\in W^{1,\infty}(\Ri^d)$ one has bounds $|c_{ij}(x)-c_{ij}(y)|\leq a\,(|x-y|\wedge 1)$ for some $a>0$.
Therefore it follows that $Q_\kappa$  
 satisfy bounds $\|Q_\kappa\|_{2\to2}\leq b\,\kappa^{-1/2}$ for all $\kappa\geq 1$.
Thus $\|Q_\kappa\|_{2\to2}<1$, the operator $I+Q_\kappa$ has a bounded inverse and 
\[
(\kappa I+\overline L)R_\kappa(I+Q_\kappa)^{-1}=I
\]
for all large $\kappa$.
Then the range of $(\kappa I+\overline L)$ is $L_2(\Ri^d)$ and $\overline L$ is self-adjoint.

\smallskip

Secondly, to deduce that $D(\overline L)=W^{2,2}(\Ri^d)$, with equivalent norms,
we note that 
\begin{eqnarray*}
\|L\varphi\|_2^2=\sum^d_{i,j,k,l=1}(\partial_j\,c_{ij}\,\partial_i\varphi, \partial_k\,c_{kl}\,\partial_l\varphi)
=\sum^d_{i,j,k,l=1}(\partial_k\partial_i\varphi, c_{ij}\,c_{kl}\,\partial_j\partial_l\varphi)+\mbox{LOT}
\end{eqnarray*}
for all $\varphi\in C_c^\infty(\Ri^d)$ where  LOT denotes  a sum of lower order terms.
But if $A$ denotes  the $d^2\times d^2$-matrix with coefficients $a_{(ik),(jl)}= c_{ij}c_{kl}$ then  $A=C\otimes C$.
Thus if $\lambda \,I\geq C\geq \mu \,I>0$ then   $\lambda^2 \,I\geq A\geq \mu^2 \,I$ uniformly on $\Ri^d$ and 
\begin{equation}
\lambda^2\,\|\Delta\varphi\|_2^2\geq \sum^d_{i,j,k,l=1}(\partial_k\partial_i\varphi, c_{ij}\,c_{kl}\,\partial_j\partial_l \varphi)
\geq \mu^2\,\|\Delta\varphi\|_2^2
\label{ese1.2}
\end{equation}
for all $\varphi\in C_c^\infty(\Ri^d)$ where $\Delta$ denotes the usual Laplacian. 
The lower order terms can, however, be bounded  by standard estimates.
For each $\varepsilon\in\langle0,1]$ there is a $c_\varepsilon>0$ such that  
\begin{equation}
|\mbox{LOT}|\leq \varepsilon\,\|\Delta\varphi\|_2^2+c_\varepsilon\,\|\varphi\|_2^2
\label{ese1.3}
\end{equation}
for all $\varphi\in C_c^\infty(\Ri^d)$.
The identity $D(\overline L)=W^{2,2}(\Ri^d)$ follows straightforwardly by combination of the estimates
 (\ref {ese1.2}) and (\ref {ese1.3})
since  the $W^{2,2}$-norm and the graph norm on $D(\Delta)$ are equivalent and $C_c^\infty(\Ri^d)$ is dense in $W^{2,2}(\Ri^d)$.

\smallskip

Thirdly,   if $\delta\geq0$ then $W^{\delta,2}(\Ri^d)=D((I+\Delta)^{\delta/2})$ with the graph norm and $W^{-\delta,2}(\Ri^d)$
is defined as the dual space. 
These spaces are a scale of spaces for real interpolation.
Since, by the foregoing, $(I+\overline L)^{-1}L_2(\Ri^d)=D(\overline L)=W^{2,2}(\Ri^d)$, 
i.e.\  the resolvent is a bounded operator from $L_2(\Ri^d)$ to $W^{2,2}(\Ri^d)$,  
 it suffices to prove that 
 $(I+\overline L)^{-1}$  is bounded from
$W^{-1,2}(\Ri^d)$ to $W^{1,2}(\Ri^d)$.
Therefore it suffices to prove that $(I+\Delta)^{1/2}(I+\overline L)^{-1}(I+\Delta)^{1/2}$ extends to a  bounded operator on $L_2(\Ri^d)$.
This follows, however, because 
\begin{eqnarray*}
\|(I+\Delta)^{1/2}(I+\overline L)^{-1}(I+\Delta)^{1/2}\varphi\|_2^2&\leq&(1\vee\mu^{-1})\,(\varphi,(I+\Delta)^{1/2}(I+\overline L)^{-1}(I+\Delta)^{1/2}\varphi)\\
&\leq &(1\vee\mu^{-1})\,\|\varphi\|_2\,\|(I+\Delta)^{1/2}(I+\overline L)^{-1}(I+\Delta)^{1/2}\varphi\|_2
\end{eqnarray*}
 by strong ellipticity.
\hfill$\Box$

\newpage


\begin{thebibliography}{ERSZ07}

\bibitem[AH05]{AH}
{\sc Ariyoshi, T., {\rm and} Hino, M.}, Small-time asymptotic estimates in
  local Dirichlet spaces.
\newblock {\em Elec.\ J. Prob.} {\bf 10} (2005),  1236--1259.

\bibitem[AKR90]{AKR}
{\sc Albeverio, S., Kusuoka, S., {\rm and} R{\"o}ckner, M.}, On partial
  integration in infinite-dimensional space and applications to Dirichlet
  forms.
\newblock {\em Proc.\ London Math.\ Soc.} {\bf 42} (1990),  122--136.

\bibitem[Aze74]{Aze}
{\sc Azencott, R.}, Behavior of diffusion semi-groups at infinity.
\newblock {\em Bull.\ Soc.\ Math.\ France} {\bf 102} (1974),  193--240.

\bibitem[BH91]{BH}
{\sc Bouleau, N., {\rm and} Hirsch, F.}, {\em Dirichlet forms and analysis on
  Wiener space}, vol.\ 14 of de Gruyter Studies in Mathematics.
\newblock Walter de Gruyter \& Co., Berlin, 1991.

\bibitem[BM95]{BM}
{\sc Biroli, M., {\rm and} Mosco, U.}, A Saint-Venant type principle for
  Dirichlet forms on discontinuous media.
\newblock {\em Ann.\ Mat.\ Pura Appl.} {\bf 169} (1995),  125--181.

\bibitem[CMP98]{CMP}
{\sc Campiti, M., Metafune, G., {\rm and} Pallara, D.}, Degenerate self-adjoint
  evolution equations on the unit interval.
\newblock {\em Semigroup Forum} {\bf 57} (1998),  1--36.

\bibitem[Dav85]{Dav14}
{\sc Davies, E.~B.}, $L^1$ properties of second order elliptic operators.
\newblock {\em Bull.\ London Math.\ Soc.} {\bf 17} (1985),  417--436.

\bibitem[Dav92]{Dav12}
\leavevmode\vrule height 2pt depth -1.6pt width 23pt, Heat kernel bounds,
  conservation of probability and the Feller property.
\newblock {\em J. Anal.\ Math.} {\bf 58} (1992),  99--119.
\newblock Festschrift on the occasion of the 70th birthday of Shmuel Agmon.

\bibitem[Ebe99]{Ebe}
{\sc Eberle, A.}, {\em Uniqueness and non-uniqueness of semigroups generated by
  singular diffusion operators}.
\newblock Lect.\ Notes in Math. 1718. Springer-Verlag, Berlin etc., 1999.

\bibitem[ER09a]{ER30}
{\sc Elst, A. F.~M. ter, {\rm and} Robinson, D.~W.}, Conservation and
  invariance properties of submarkovian semigroups.
\newblock {\em J. Ramanujan.\ Math.\ Soc.} {\bf 24} (2009),  1--13.

\bibitem[ER09b]{ER27}
\leavevmode\vrule height 2pt depth -1.6pt width 23pt, Uniform subellipticity.
\newblock {\em J. Operator Theory} {\bf 62} (2009),  125--149.

\bibitem[ERSZ06]{ERSZ2}
{\sc Elst, A. F.~M. ter, Robinson, D.~W., Sikora, A., {\rm and} Zhu, Y.},
  Dirichlet forms and degenerate elliptic operators.
\newblock In {\sc Koelink, E., Neerven, J. van, Pagter, B. de, {\rm and}
  Sweers, G.}, eds., {\em Partial Differential Equations and Functional
  Analysis}, vol.\ 168 of Operator Theory: Advances and Applications.
  Birkh{\"a}user, 2006,  73--95.
\newblock Philippe Clement Festschrift.

\bibitem[ERSZ07]{ERSZ1}
{\sc Elst, A. F.~M. ter, Robinson, D.~W., Sikora, A., {\rm and} Zhu, Y.},
  Second-order operators with degenerate coefficients.
\newblock {\em Proc.\ London Math.\ Soc.} {\bf 95} (2007),  299--328.

\bibitem[Fal03]{Fal}
{\sc Falconer, K.~J.}, {\em Fractal geometry}.
\newblock Second edition, Mathematical foundations and applications. John Wiley
  \& Sons Inc., Hoboken, NJ, 2003.

\bibitem[FOT94]{FOT}
{\sc Fukushima, M., Oshima, Y., {\rm and} Takeda, M.}, {\em Dirichlet forms and
  symmetric Markov processes}, vol.\ 19 of de Gruyter Studies in Mathematics.
\newblock Walter de Gruyter \& Co., Berlin, 1994.

\bibitem[Fri34]{Friedr2}
{\sc Friedrichs, K.~O.}, Spektraltheorie halbbeschr{\"a}nkter Operatoren und
  Anwendung auf die Spektralzerlegung von Differentialoperatoren. I.
\newblock {\em Math.\ Anal.} {\bf 109} (1934),  465--487.

\bibitem[Fri64]{Fri}
{\sc Friedman, A.}, {\em Partial differential equations of parabolic type}.
\newblock Prentice-Hall, Inc., Englewood Cliffs, N.J., 1964.

\bibitem[Gaf59]{Gaf}
{\sc Gaffney, M.~P.}, The conservation property of the heat equation on
  Riemannian manifolds.
\newblock {\em Comm.\ Pure Appl.\ Math.} {\bf 12} (1959),  1--11.

\bibitem[GT83]{GT}
{\sc Gilbarg, D., {\rm and} Trudinger, N.~S.}, {\em Elliptic partial
  differential equations of second order}.
\newblock Second edition, Grundlehren der mathematischen Wissenschaften 224.
  Springer-Verlag, Berlin etc., 1983.

\bibitem[Kat80]{Kat1}
{\sc Kato, T.}, {\em Perturbation theory for linear operators}.
\newblock Second edition, Grundlehren der mathematischen Wissenschaften 132.
  Springer-Verlag, Berlin etc., 1980.

\bibitem[Kre47]{Kre1}
{\sc Krein, M.~G.}, The theory of self-adjoint extensions of semi-bounded
  Hermitian transformations and its applications. I.
\newblock {\em Mat.\ Sbornik N.S.} {\bf 20(62)} (1947),  431--495.

\bibitem[MR92]{MR}
{\sc Ma, Z.~M., {\rm and} R{\"o}ckner, M.}, {\em Introduction to the theory of
  (non symmetric) Dirichlet Forms}.
\newblock Universitext. Springer-Verlag, Berlin etc., 1992.

\bibitem[Ouh05]{Ouh5}
{\sc Ouhabaz, E.-M.}, {\em Analysis of heat equations on domains}, vol.\ 31 of
  London Mathematical Society Monographs Series.
\newblock Princeton University Press, Princeton, NJ, 2005.

\bibitem[Rob87]{Rob7}
{\sc Robinson, D.~W.}, Commutator theory on Hilbert space.
\newblock {\em Can.\ J. Math.} {\bf 34} (1987),  1235--1280.

\bibitem[RS07]{RSi}
{\sc Robinson, D.~W., {\rm and} Sikora, A.}, Degenerate elliptic operators:
  capacity, flux and separation.
\newblock {\em J. Ramanujan Math.\ Soc.} {\bf 22} (2007),  385--408.

\bibitem[RS08]{RSi2}
\leavevmode\vrule height 2pt depth -1.6pt width 23pt, Analysis of degenerate
  elliptic operators of Gru\v{s}in type.
\newblock {\em Math.\ Z.} {\bf 260} (2008),  475--508.

\bibitem[RS09]{RSi3}
\leavevmode\vrule height 2pt depth -1.6pt width 23pt, Degenerate elliptic
  operators in one-dimension.
\newblock Research report, Australian National University, 2009.
\newblock arXiv:math./AP0909.0567.

\bibitem[Stu98]{Stu2}
{\sc Sturm, K.-T.}, The geometric aspect of Dirichlet forms.
\newblock In {\em New directions in Dirichlet forms}, vol.\ 8 of AMS/IP Stud.\
  Adv.\ Math.,  233--277. Amer.\ Math.\ Soc., Providence, RI, 1998.

\bibitem[Tak96]{Take}
{\sc Takeda, M.}, Two classes of extensions for generalized {S}chr{\"o}dinger
  operators.
\newblock {\em Pot.\ Anal.} {\bf 5} (1996),  1--13.

\end{thebibliography}
\end{document}